\documentclass[reqno]{amsart}
\usepackage[T1]{fontenc}
\usepackage{mathptmx}
\usepackage{amssymb}
\usepackage{color,graphics}
\usepackage[colorlinks]{hyperref}
\usepackage[all]{xy}

\newcommand{\CC}{\mathcal{C}}

\newcommand{\BB}{\mathcal{B}}
\newcommand{\EE}{\mathcal{E}}

\newcommand{\DD}{\mathcal{D}}
\newcommand{\UU}{\mathcal{U}}
\newcommand{\VV}{\mathcal{V}}
\newcommand{\WW}{\mathcal{W}}

\newcommand{\OO}{\mathcal{O}}
\renewcommand{\SS}{\mathcal{S}}
\newcommand{\RR}{\mathcal{R}}
\newcommand{\topcat}{\textnormal{\textbf{Top}}}
\newcommand{\diffcat}{\textnormal{\textbf{Diff}}}

\newcommand{\mypath}{\textnormal{\textbf{Path}}}

\newcommand{\myop}{\textnormal{op}}
\renewcommand{\hom}{\textnormal{Hom}}
\newcommand{\myequiv}{\simeq}

\newcommand{\limit}{\varprojlim}

\newcommand{\myemph}[1]{\textnormal{\textbf{#1}}}
\newcommand{\sheaves}{\textnormal{\textbf{Sh}}}
\newcommand{\stacks}{\textnormal{\textbf{St}}}
\newcommand{\cat}{\textnormal{\textbf{Cat}}}
\newcommand{\sets}{\textnormal{\textbf{Set}}}
\newcommand{\iso}{\cong}
\newcommand{\pseudo}[2]{{[}#1,#2{]}_{\textnormal{ps}}}
\newcommand{\prestacks}{\textnormal{\textbf{PreSt}}}
\newcommand{\toppre}{\textnormal{\textbf{TopPreSt}}}
\newcommand{\topstacks}{\textnormal{\textbf{TopSt}}}
\newcommand{\diffpre}{\textnormal{\textbf{DiffPreSt}}}
\newcommand{\diffstacks}{\textnormal{\textbf{DiffSt}}}
\newcommand{\schemes}{\textnormal{\textbf{Sch}}}
\newcommand{\algpre}{\textnormal{\textbf{AlgPreSt}}}
\newcommand{\algstacks}{\textnormal{\textbf{AlgSt}}}
\newcommand{\hocat}{\textnormal{\textbf{Ho}}}
\newcommand{\descent}[2]{\textnormal{\textbf{Desc}}(#1,#2)}
\newcommand{\astack}{\textnormal{\textbf{a}}}

\newcommand{\groupoids}{\textnormal{\textbf{Gpd}}}

\newtheorem{theorem}{Theorem}[section]
\newtheorem{lemma}[theorem]{Lemma}
\newtheorem{proposition}[theorem]{Proposition}

\newtheorem*{corollary}{Corollary}
\newtheorem*{theorem*}{Theorem}

\theoremstyle{definition}
\newtheorem{definition}{Definition}[section]
\newtheorem{example}[definition]{Example}

\theoremstyle{remark}
\newtheorem*{remark}{Remark}

\title{bicategorical fibration structures and stacks}
\author{Dorette A. Pronk and Michael A. Warren}
\address{Department of Mathematics and Statistics\\Dalhousie
  University\\Halifax, Nova Scotia, Canada}
\address{School of Mathematics\\Institute for Advanced
  Study\\Princeton, New Jersey, USA}
\date{\today}
\keywords{stacks, fibrant objects, homotopy bicategory,
  bicategories of fractions, algebraic stacks, differentiable stacks,
  topological stacks}
\subjclass[2010]{Primary: 18D05; Secondary: 18G55, 14A20}

\begin{document}

\maketitle

\begin{abstract}
  The familiar construction of categories of fractions, due to Gabriel
  and Zisman, allows one to invert a class $W$ of arrows in a category
  in a universal way. Similarly, bicategories of fractions allow one
  to invert a collection $W$ of arrows in a bicategory $\BB$. In this
  case the arrows are inverted in the sense that they are made into
  equivalences. As with categories of fractions, bicategories of
  fractions suffer from the defect that they need not be locally small
  even when $\BB$ is locally small. Similarly, in the case where $\BB$
  is a 2-category, the bicategory of fractions will not in general be
  a 2-category.

  In this paper we introduce two notions ---systems of fibrant objects and
  fibration systems--- which will allow us to associate to a bicategory $\BB$ a
  homotopy bicategory $\hocat(\BB)$ in such a way that $\hocat(\BB)$
  is the universal way to invert weak equivalences in $\BB$. This
  construction resolves both of the difficulties with bicategories of
  fractions mentioned above. We also describe a fibration system on
  the 2-category of prestacks on a site and prove that the resulting
  homotopy bicategory is the 2-category of stacks. Further examples
  considered include algebraic, differentiable and topological
  stacks.
\end{abstract}

\section*{Introduction}

It is widely known that Quillen's \cite{Quillen:HA} notion of
\emph{model structure} provides an adequate (for many purposes) formal
setting for the development of the theory of $(\infty,1)$-categories,
as has been studied by Joyal \cite{Joyal:TQCA}, Lurie \cite{Lurie:HTT}
and others.  Moreover, a model structure on a category $\CC$ provides
a technical tool for forming the localization of $\CC$ with respect to
a class of weak equivalences: weak equivalences are inverted in a
universal way in the passage to the homotopy category $\hocat(\CC)$ of
$\CC$.  Consequently, it is possible to invert weak equivalences in
this setting without having to resort to the Gabriel-Zisman
\cite{Gabriel:CFHT} calculus of fractions.

In the bicategorical setting, one might like to be able to invert a collection 
of weak equivalences in the sense of turning them into
equivalences.  In
\cite{Pronk:PhD,Pronk:ESBF}, the first author gave a bicategorical
generalization of the Gabriel-Zisman calculus of fractions which
accomplishes this goal:
\begin{theorem*}[Pronk \cite{Pronk:ESBF}]
  Given a collection of arrows $\frak{W}$ in a bicategory $\CC$
  satisfying certain conditions, there exists a bicategory
  $\CC[\frak{W}^{-1}]$ (called the \emph{bicategory of fractions for $\frak{W}$}) and
  a homomorphism $I\colon \CC\to\CC[\frak{W}^{-1}]$ such that $I$ sends
  arrows in $\frak{W}$ to equivalences in $\CC[\frak{W}^{-1}]$ and $I$ is universal
  with this property.
\end{theorem*}

Like the ordinary category of fractions, this construction suffers
from the technical defect that $\CC[\frak{W}^{-1}]$ will not in
general have small hom-categories even when $\CC$ does.  Moreover,
$\CC[\frak{W}^{-1}]$ will be a bicategory even when $\CC$ is a
2-category.  In the present paper, we introduce the notion of a
\emph{system of fibrant objects} in a bicategory $\CC$ and the notion
of a \emph{fibration structure} on a bicategory which will allow us to
form the localization of a bicategory $\CC$ with respect to a class of weak
equivalences in such a way that the result will both have small
hom-categories when $\CC$ does and will be a 2-category when $\CC$ is.

A system of fibrant objects consists of a collection $\frak{W}$
(\emph{weak equivalences}) of maps in $\CC$, a pseudofunctor
$Q\colon\CC\to\CC$ (fibrant replacement) and a pseudonatural transformation $\eta\colon
1_{\CC}\to Q$ (whose components are weak equivalences) satisfying certain factorization
conditions.  To each bicategory $\CC$ with a
system of fibrant objects, there is an
associated bicategory $\hocat(\CC)$, called the \emph{homotopy
  bicategory of $\CC$}, and a pseudofunctor
$I\colon\CC\to\hocat(\CC)$.  By construction, $\hocat(\CC)$ has small
hom-categories when $\CC$ does and it is a 2-category when $\CC$
is. Our first main result is as follows:
\begin{theorem*}[Theorem \ref{theorem:main} below]
  The pseudofunctor $I$ inverts weak equivalences and is universal
  with this property.
\end{theorem*}

The remainder of this paper is concerned with investigating specific
examples of bicategories with systems of fibrant objects.  Our leading
example is the 2-category $\stacks(\CC)$ of stacks on a site $\CC$,
which we obtain (Corollary \ref{theorem:stacks}) as the homotopy
2-category $\hocat(\prestacks(\CC))$ of the 2-category of prestacks on
$\CC$.  This result is made possible using a characterization of the
fibrations of prestacks which is analogous to an earlier
characterization of stacks given, albeit in a different setting from the one
considered here, by Joyal and Tierney \cite{Joyal:SSCS} (cf.~also the
work of Hollander \cite{Hollander}).

The system of fibrant objects on $\prestacks(\CC)$ is notable in that
it exhibits a number of additional features making it more closely
resemble the notion of a model structure.  These additional features
are sufficiently interesting that we introduce the notion of a
\emph{fibration structure} on a bicategory to capture them.

A category $\CC$ has a fibration structure when there are stronger lifting
and factorization conditions in place which among other things imply that
the category has path objects and that the factorization lemma holds, so 
that one can construct generalized universal bundles.

In \cite{Pronk:ESBF}, the first author gave a number of
examples of bicategorical equivalences between well-known 2-categories
and bicategories of fractions.  These examples include topological,
differentiable and algebraic stacks and we show that these examples
can also be captured in our setting.  Note however that the characterizations
given here of these 2-categories differ from those in \emph{ibid}.
In \emph{ibid} these 2-categories were characterized as bicategories of fractions
of certain categories of groupoids with respect to Morita equivalences.
Here we will view them as homotopy categories of certain categories of prestacks
with respect to local weak equivalences.

Part of the motivation of this work is the goal of trying to find a
formal setting, analogous to the setting of model structures, in which
to develop the theory of $(\infty,2)$- and $(\infty,n)$-categories.
Ultimately we would like to extend the axiomatization given here to
the lax setting (we are always working in a ``pseudo'' setting) and to
relate the results presented here to Street's notion of 2-topos
\cite{Street:TDST}.  Intuitively, every 2-topos should arise as a
homotopy 2-category by analogy with the way Grothendieck toposes arise as
localizations of presheaf categories.

\subsection*{Summary}

In Section \ref{section:basics} we recall basic definitions and
results on bicategories, pseudofunctors, pseudonatural
transformations, and so forth.  In Section
\ref{section:fibration_structures} we introduce systems of fibrant
objects and fibration structures on bicategories and we prove our main
result (Theorem \ref{theorem:main}).  In Section \ref{section:local_fib} we introduce a fibered notion of stack:
\emph{local fibrations}.  Let $(\CC,J)$ be a site and let
pseudofunctors $E,B\colon \CC^{\myop}\to\cat$ and a pseudonatural
transformation $p\colon E\to B$ be given.  For each cover $\SS$ of an object
$U$ of $\CC$
we introduce the category $\descent{p}{\SS}$ of \emph{descent data
  with respect to $p$ and $\SS$}.  This category, like the usual
category of descent data $\descent{E}{\SS}$, can be defined as a
pseudolimit (although here we give a direct description) and we define $p$ to be a local fibration when it
satisfies an effective descent condition with respect to
$\descent{p}{\SS}$ analogous to the usual descent condition for
stacks.  In Section \ref{section:model_structure} we describe a
fibration structure on the 2-category of prestacks $\prestacks(\CC)$
and prove that the resulting homotopy bicategory is the 2-category
$\stacks(\CC)$ of stacks.  In particular, we introduce the local weak
equivalences (which are already known in the literature on stacks) and
prove, using the Axiom of Choice, that the local fibrations are exactly those maps having a
bicategorical version of the right lifting property with respect to
the local weak equivalences.  Further examples (algebraic,
differentiable and topological prestacks) of systems of fibrant
objects are considered in Section \ref{section:examples}.

\section{Basics and notation}\label{section:basics}

It is worth mentioning that we make free use of the Axiom of Choice.
As such, we do not distinguish between \emph{strong} and \emph{weak
  categorical equivalences} (see \cite{Bunge:SEIC} for more on the
differences between strong and weak equivalences).   We assume that the
reader is familiar with the basic theory of 2-categories and refer the
reader to \cite{Kelly:RE2C} for further details. For more information
regarding stacks we refer the reader to \cite{Giraud:CNA} and \cite{Moerdijk:ILSG}.  

\subsection{Bicategories}\label{section:bicats}

We now review the definitions of bicategories, pseudofunctors,
pseudonatural transformations and modifications.
\begin{definition}[B\'{e}nabou \cite{Benabou:IB}]\label{def:bicat}
  A \myemph{bicategory} $\CC$ consists of a collection of objects
  $A,B,\ldots$ together with the following data:
  \begin{itemize}
    \item Categories $\CC(A,B)$ for objects $A$ and $B$ of $\CC$.  The
      objects of $\CC(A,B)$ are called \emph{arrows} and the
      arrows are called \emph{2-cells}.  When $\alpha$ and $\beta$ are
      composable 2-cells in $\CC(A,B)$ we denote their composite by
      $\beta\cdot\alpha$.
    \item For objects $A,B$ and $C$ of $\CC$, a functor
      $c_{A,B,C}:\CC(A,B)\times\CC(B,C)\to\CC(A,C)$.  We denote
      $c_{A,B,C}(f,g)$ by $g\circ f$, for arrows $f\in\CC(A,B)$ and
      $g\in\CC(B,C)$, and we denote $c_{A,B,C}(\alpha,\beta)$ by
      $\beta*\alpha$, for 2-cells $\beta\in\CC(A,B)$ and
      $\alpha\in\CC(B,C)$.  When no confusion will result we omit the
      subscripts and write $c$ instead of $c_{A,B,C}$.
    \item For each object $A$ of $\CC$, an arrow $1_{A}\in\CC(A,A)$.
    \item For objects $A,B,C$ and $D$ of $\CC$, a natural isomorphism:
      \begin{align*}
        \xy
        {\ar^(.55){c\times\CC(C,D)}(0,15)*+{\CC(A,B)\times\CC(B,C)\times\CC(C,D)};(50,15)*+{\CC(A,C)\times\CC(C,D)}};
        {\ar_{\CC(A,B)\times
            c}(0,15)*+{\CC(A,B)\times\CC(B,C)\times\CC(C,D)};(0,0)*+{\CC(A,B)\times\CC(B,D)}};
        {\ar_{c}(0,0)*+{\CC(A,B)\times\CC(B,D)};(50,0)*+{\CC(A,D).}};
        {\ar^{c}(50,15)*+{\CC(A,C)\times\CC(C,D)};(50,0)*+{\CC(A,D).}};
        {\ar@{=>}^{\alpha_{A,B,C,D}}(22,7.5)*+{};(28,7.5)*+{}};
        \endxy
      \end{align*} 
      As with the composition functors $c$, we will omit subscripts
      and write $\alpha$ instead of $\alpha_{A,B,C,D}$.
    \item For objects $A$ and $B$ of $\CC$, natural isomorphisms
      $\lambda_{A,B}$ and $\rho_{A,B}$ as
      indicated in the following diagrams:
      \begin{align*}
        \xy
        {\ar^(.4){1_{A}\times\CC(A,B)}(0,15)*+{1\times\CC(A,B)};(36,15)*+{\CC(A,A)\times\CC(A,B)}};
        {\ar_{\pi_{1}}(0,15)*+{1\times\CC(A,B)};(18,0)*+{\CC(A,B)}};
        {\ar^{c}(36,15)*+{\CC(A,A)\times\CC(A,B)};(18,0)*+{\CC(A,B)}};
        {\ar@{<=}^{\rho_{A,B}}(15,7.5)*+{};(21,7.5)*+{}};
        \endxy
        \quad
        \xy
        {\ar^(.4){\CC(A,B)\times 1_{B}}(0,15)*+{\CC(A,B)\times 1};(36,15)*+{\CC(A,B)\times\CC(B,B)}};
        {\ar_{\pi_{0}}(0,15)*+{\CC(A,B)\times 1};(18,0)*+{\CC(A,B).}};
        {\ar^{c}(36,15)*+{\CC(A,B)\times\CC(B,B)};(18,0)*+{\CC(A,B).}};
        {\ar@{<=}^{\lambda_{A,B}}(15,7.5)*+{};(21,7.5)*+{}};
        \endxy
      \end{align*}
      We again omit subscripts and simply write $\lambda$ and $\rho$.
    \end{itemize}
    These data are required to satisfy the following conditions:
    \begin{itemize}
    \item Given $f\in\CC(A,B)$, $g\in\CC(B,C)$, $h\in\CC(C,D)$ and
      $k\in\CC(D,E)$, the following diagram commutes:
      \begin{align*}
        \xy
        {\ar^{\alpha*f}(0,15)*+{\bigl((k\circ h)\circ g\bigr)\circ
            f};(60,15)*+{\bigl(k\circ (h\circ g)\bigr)\circ f}};
        {\ar_{\alpha}(0,15)*+{\bigl((k\circ h)\circ g\bigr)\circ
            f};(0,0)*+{(k\circ h)\circ (g\circ f)}};
        {\ar_{\alpha}(0,0)*+{(k\circ h)\circ (g\circ
            f)};(30,0)*+{k\circ\bigl(h\circ(g\circ f)\bigr)}};
        {\ar^{\alpha}(60,15)*+{\bigl(k\circ (h\circ g)\bigr)\circ
            f};(60,0)*+{k\circ\bigl((h\circ g)\circ f\bigr).}};
        {\ar^{k*\alpha}(60,0)*+{k\circ\bigl((h\circ g)\circ
            f\bigr).};(30,0)*+{k\circ\bigl(h\circ (g\circ f)\bigr)}};
        \endxy
      \end{align*}
    \item Given $f\in\CC(A,B)$ and $g\in\CC(B,C)$, the following
      diagram commutes:
      \begin{align*}
        \xy
        {\ar^{\alpha}(0,15)*+{(g\circ 1_{B})\circ
            f};(30,15)*+{g\circ(1_{B}\circ f)}};
        {\ar_{\rho*f}(0,15)*+{(g\circ 1_{B})\circ f};(15,0)*+{g\circ
            f.}};
        {\ar^{g*\lambda}(30,15)*+{g\circ(1_{B}\circ
            f)};(15,0)*+{g\circ f.}};
        \endxy
      \end{align*}
    \end{itemize}
\end{definition}
The following definition also involves coherence data which should
technically carry subscripts.  These are indicated explicitly the
first time they appear, but afterwards we adopt a policy of omitting
subscripts wherever possible as in Definition \ref{def:bicat}.
\begin{definition}
  Given bicategories $\CC$ and $\DD$, a \myemph{pseudofunctor}
  $F\colon\CC\to\DD$ is given by the following data:
  \begin{itemize}
  \item An assignment of objects $FC$ of $\DD$ to each object $C$ of $\CC$.
  \item For all objects $A$ and $B$ of $\CC$, a functor
    $F_{A,B}\colon\CC(A,B)\to\CC(FA,FB)$.
  \item For all objects $A,B$ and $C$ of $\CC$, natural isomorphisms
    as indicated in the following diagrams:
    \begin{align*}
      \xy
      {\ar^(.6){c}(0,15)*+{\CC(A,B)\times\CC(B,C)};(40,15)*+{\CC(A,C)}};
      {\ar_{F\times
          F}(0,15)*+{\CC(A,B)\times\CC(B,C)};(0,0)*+{\DD(FA,FB)\times\DD(FB,FC)}};
      {\ar_(.6){c}(0,0)*+{\DD(FA,FB)\times\DD(FB,FC)};(40,0)*+{\DD(FA,FC)}};
      {\ar^{F}(40,15)*+{\CC(A,C)};(40,0)*+{\DD(FA,FC)}};
      {\ar@{=>}^{\varphi_{A,B,C}}(17,7.5)*+{};(23,7.5)*+{}};
      \endxy
    \end{align*}
    and
    \begin{align*}
      \xy
      {\ar^{1_{A}}(0,15)*+{1};(30,15)*+{\CC(A,A)}};
      {\ar_{1_{FA}}(0,15)*+{1};(15,0)*+{\DD(FA,FA).}};
      {\ar^{F}(30,15)*+{\CC(A,A)};(15,0)*+{\DD(FA,FA).}};
      {\ar@{=>}^(.35){\upsilon_{A}}(12,7.5)*+{};(18,7.5)*+{}};
      \endxy
    \end{align*}
  \end{itemize}
  These data are required to be such that the following diagrams
  commute:
  \begin{align*}
    \xy
    {\ar^{\varphi*Ff}(0,15)*+{(Fh\circ Fg)\circ Ff};(30,15)*+{F(h\circ
        g)\circ Ff}};
    {\ar^{\varphi}(30,15)*+{F(h\circ g)\circ
        Ff};(60,15)*+{F\bigl((h\circ g)\circ f\bigr)}};
    {\ar_{\alpha}(0,15)*+{(Fh\circ Fg)\circ
        Ff};(0,0)*+{Fh\circ(Fg\circ Ff)}};
    {\ar_{Fh*\varphi}(0,0)*+{Fh\circ(Fg\circ Ff)};(30,0)*+{Fh\circ
        F(g\circ f)}};
    {\ar_{\varphi}(30,0)*+{Fh\circ F(g\circ
        f)};(60,0)*+{F\bigl(h\circ(g\circ f)\bigr)}};
    {\ar^{F\alpha}(60,15)*+{F\bigl((h\circ g)\circ f\bigr)};(60,0)*+{F\bigl(h\circ(g\circ f)\bigr)}};
    \endxy
  \end{align*}
  \begin{align*}
    \xy
    {\ar^(.45){Ff*\upsilon}(0,15)*+{Ff\circ 1_{FA}};(20,15)*+{Ff\circ
        F1_{A}}};
    {\ar^(.45){\varphi}(20,15)*+{Ff\circ F1_{A}};(40,15)*+{F(f\circ
        1_{A})}};
    {\ar_{\rho}(0,15)*+{Ff\circ 1_{FA}};(20,0)*+{Ff}};
    {\ar^{F\rho}(40,15)*+{F(f\circ 1_{A})};(20,0)*+{Ff}};
    \endxy
    \quad\text{and}\quad
    \xy
    {\ar^{\upsilon*Ff}(0,15)*+{1_{FB}\circ
        Ff};(20,15)*+{F1_{B}\circ Ff}};
    {\ar^(.45){\varphi}(20,15)*+{F1_{B}\circ
        Ff};(40,15)*+{F(1_{B}\circ f)}};
    {\ar_{\lambda}(0,15)*+{1_{FB}\circ Ff};(20,0)*+{Ff}};
    {\ar^{F\lambda}(40,15)*+{F(1_{B}\circ f)};(20,0)*+{Ff}};
    \endxy
  \end{align*}
  for any arrows $f\colon A\to B$, $g\colon B\to C$ and $h\colon C\to
  D$ in $\CC$.
\end{definition}
\begin{definition}
  Given pseudofunctors $F,G\colon\CC\to\DD$, a \myemph{pseudonatural
    transformation} $\psi\colon F\to G$ consists of the following
  data:
  \begin{itemize}
  \item For each object $C$ of $\CC$, an arrow $\psi_{C}\colon FC\to
    GC$ in $\DD$.
  \item For objects $A$ and $B$ of $\CC$, a natural isomorphism 
    \begin{align*}
      \xy
      {\ar^{F}(0,15)*+{\CC(A,B)};(30,15)*+{\DD(FA,FB)}};
      {\ar_{G}(0,15)*+{\CC(A,B)};(0,0)*+{\DD(GA,GB)}};
      {\ar^{\DD(FA,\psi)}(30,15)*+{\DD(FA,FB)};(30,0)*+{\DD(FA,GB).}};
      {\ar_{\DD(\psi,GB)}(0,0)*+{\DD(GA,GB)};(30,0)*+{\DD(FA,GB).}};
      {\ar@{=>}^{\psi_{A,B}}(12,7.5)*+{};(18,7.5)*+{}};
      \endxy
    \end{align*}
    Note that here we are overloading the notation $\psi$.  In
    practice this should not result in any confusion.
  \end{itemize}
  These data are required to be such that the following diagrams commute:
  \begin{align*}
    \xy
    {\ar^{\alpha}(0,15)*+{(Gg\circ
        Gf)\circ\psi};(30,15)*+{Gg\circ(Gf\circ\psi)}};
    {\ar^{Gg*\psi}(30,15)*+{Gg\circ(Gf\circ\psi)};(60,15)*+{Gg\circ(\psi\circ
        Ff)}};
    {\ar^{\alpha^{-1}}(60,15)*+{Gg\circ(\psi\circ
        Ff)};(90,15)*+{(Gg\circ\psi)\circ Ff}};
    {\ar^{\psi*Ff}(90,15)*+{(Gg\circ\psi)\circ Ff};(90,0)*+{(\psi\circ
        Fg)\circ Ff}};
    {\ar_{\alpha}(90,0)*+{(\psi\circ Fg)\circ Ff};(60,0)*+{\psi\circ
        (Fg\circ Ff)}};
    {\ar_{\psi*\varphi}(60,0)*+{\psi\circ (Fg\circ
        Ff)};(30,0)*+{\psi\circ F(g\circ f)}};
    {\ar_{\varphi*\psi}(0,15)*+{(Gg\circ
        Gf)\circ\psi};(0,0)*+{G(g\circ f)\circ\psi}};
    {\ar_{\psi}(0,0)*+{G(g\circ f)\circ\psi};(30,0)*+{\psi\circ
        F(g\circ f)}};
    \endxy
  \end{align*}
  \begin{align*}
    \xy
    {\ar^{\lambda}(0,15)*+{1_{GA}\circ\psi};(30,15)*+{\psi}};
    {\ar^{\rho^{-1}}(30,15)*+{\psi};(60,15)*+{\psi\circ 1_{FA}}};
    {\ar^{\psi*\upsilon}(60,15)*+{\psi\circ 1_{FA}};(60,0)*+{\psi\circ
        F1_{A}.}};
    {\ar_{\upsilon*\psi}(0,15)*+{1_{GA}\circ\psi};(0,0)*+{G1_{A}\circ\psi}};
    {\ar_{\psi}(0,0)*+{G1_{A}\circ\psi};(60,0)*+{\psi\circ F1_{A}.}};
    \endxy
  \end{align*}
\end{definition}

\begin{definition}
  A \myemph{modification} $\omega\colon\psi\to\psi'$, for $\psi$ and
  $\psi'$ pseudonatural transformations $F\to G$, consists of
  an assignment of 2-cells $\omega_{C}\colon \psi_{C}\to
  \psi'_{C}$ to each object $C$ of $\CC$ such that 
  \begin{align*}
    \xy
    {\ar^{Gf*\omega}(0,15)*+{Gf*\psi};(30,15)*+{Gf*\psi'}};
    {\ar_{\psi}(0,15)*+{Gf*\psi};(0,0)*+{\psi*Ff}};
    {\ar_{\omega*Ff}(0,0)*+{\psi*Ff};(30,0)*+{\psi'*Ff}};
    {\ar^{\psi'}(30,15)*+{Gf*\psi'};(30,0)*+{\psi'*Ff}};
    \endxy
  \end{align*}
  commutes, for each $f\colon A\to B$ in $\CC$.
\end{definition}

\subsection{Pseudofunctor bicategories}

Given bicategories $\CC$ and $\DD$, we denote by $[\CC,\DD]$ the
bicategory which has as objects pseudofunctors $\CC\to\DD$, as arrows
pseudonatural transformations, and as 2-cells modifications.  Note
that $[\CC,\DD]$ is a 2-category when $\DD$ is.
\begin{definition}
  An arrow $f\colon A\to B$ in a bicategory $\CC$ is an
  \myemph{equivalence} if there exists an arrow $f'\colon B\to A$
  together with invertible 2-cells $f\circ f'\iso 1_{B}$ and
  $1_{A}\iso f'\circ f$.
\end{definition}
It is a well-known fact that an equivalence of categories can always
be altered to give an adjoint equivalence.  The same fact holds in an
arbitrary bicategory:
\begin{lemma}\label{lemma:equiv2adjequiv}
  If $f\colon A\to B$ is an equivalence in a bicategory $\CC$, then
  there exists a $f'\colon B\to A$ together with invertible 2-cells $\eta\colon
  1_{A}\iso f'\circ f$ and $\epsilon\colon f\circ f' \iso 1_{B}$
  which are the unit and counit of an adjunction $f\dashv f'$.
\end{lemma}
\begin{lemma}\label{lemma:pseudonat}
  Given pseudofunctors $F,G\colon\CC\to\DD$ between bicategories $\CC$
  and $\DD$, if $\xi\colon F\to G$ is a pseudonatural transformation
  such that, for each $A$ in $\CC$, $\xi_{A}\colon FA\to GA$ is an
  equivalence, then there exists a pseudonatural transformation
  $\xi'\colon G\to F$ such that $\xi'$ is an adjoint pseudoinverse of
  $\xi$ in the bicategory $[\CC,\DD]$.
  \begin{proof}
    By Lemma \ref{lemma:equiv2adjequiv}, we may choose $\xi'_{A}$
    together with $\eta_{A}\colon 1_{FA}\iso \xi'_{A}\circ\xi_{A}$ and
    $\epsilon_{A}\colon \xi_{A}\circ\xi'_{A}\iso 1_{GA}$ making
    $\xi_{A}\dashv \xi'_{A}$.  Then, for $f\colon A\to B$ in $\CC$,
    the isomorphism $Ff\circ\xi'_{A}\iso\xi'_{B}\circ Gf$ is
    constructed by composing the isomorphisms 
    \begin{align*}
      Ff\circ\xi'_{A} \;\iso\; \xi'_{B}\circ (\xi_{B}\circ
      Ff)\circ \xi'_{A}\; \iso \; \xi'_{B}\circ (Gf\circ \xi_{A})\circ
      \xi'_{A}\;\iso\;\xi'_{B}\circ Gf
    \end{align*}
    where the first isomorphism is a result of the coherence
    isomorphisms together with $\eta_{B}$, the second isomorphism is
    by $\xi_{f}$ and the third is by coherence and $\epsilon_{A}$.
    The coherence conditions on pseudonatural transformations
    follow from pseudonaturality of $\xi$ and the triangle laws for
    adjunctions.
  \end{proof}
\end{lemma}
\begin{definition}
  Given bicategories $\CC$ and $\DD$, a pseudofunctor
  $F\colon\CC\to\DD$ is an \myemph{equivalence of bicategories} if
  there exists a pseudofunctor $G\colon\DD\to\CC$ together with maps
  $\eta\colon 1_{\CC}\to G\circ F$ and $\epsilon\colon F\circ G\to
  1_{\DD}$ which are equivalences in the bicategories $[\CC,\CC]$ and
  $[\DD,\DD]$, respectively.
\end{definition}
\begin{definition}
  A pseudofunctor $F\colon\CC\to\DD$ is a \myemph{weak equivalence of
    bicategories} if the following conditions are satisfied:
  \begin{itemize}
  \item For each object $D$ of $\DD$, there exists an object $C$ of
    $\CC$ and an equivalence $FC\to D$ in $\DD$.
  \item For all objects $C$ and $C'$ of $\CC$, the map
    $\CC(C,C')\to\DD(FC,FC')$ is an equivalence of categories.
  \end{itemize}
\end{definition}
Note that, in the presence of the Axiom of Choice, the notions of
equivalence and weak equivalence of bicategories coincide.

\subsection{Arrow bicategories}

Given a bicategory $\CC$ we define a new bicategory $\CC^{\to}$ as
follows:
\begin{description}
\item[Objects] An object is an arrow $f\colon A\to B$ in $\CC$.
\item[Arrows] Given objects $f\colon A\to B$ and $g\colon C\to D$, an
  arrow $f\to g$ is given by arrows $h\colon A\to C$ and $k\colon B\to
  D$ together with an invertible 2-cell $\gamma$ as indicated in the
  following diagram:
     \begin{align*}
      \xy
      {\ar^{h}(0,20)*+{A};(20,20)*+{C}};
      {\ar_{f}(0,20)*+{A};(0,0)*+{B}};
      {\ar^{g}(20,20)*+{C};(20,0)*+{D.}};
      {\ar_{k}(0,0)*+{B};(20,0)*+{D.}};
      {\ar^{\gamma}@{=>}(10,13)*+{};(10,7)*+{}};
      \endxy
    \end{align*}
\item[2-cells] Given objects $f$ and $g$, and arrows $(h,k,\gamma)$
  and $(h',k',\gamma')$ from $f$ to $g$, a
  2-cell $\varphi\colon (h,k,\gamma)\to(h',k',\gamma')$ consists of
  invertible 2-cells $\varphi_{0}\colon h\to h'$ and $\varphi_{1}\colon k\to k'$
  in $\CC$ such that 
  \begin{align*}
    \gamma'\cdot (g*\varphi_{0}) & = (\varphi_{1}*f)\cdot \gamma.
  \end{align*}
\item[Horizontal composition] Given objects $f\colon A\to B$, $g\colon
  C\to D$ and $i\colon E\to F$, and arrows $(h,k,\gamma)\colon f\to g$
  and $(l,m,\delta)\colon g\to i$, we define
  $(l,m,\delta)*(h,k,\gamma)$ to be the composite
$$
(l\circ h,m\circ k,\xymatrix@1{i\circ(l\circ h)\ar[r]^{\alpha^{-1}}&(i\circ l)\circ h\ar[r]^{\delta * h}
    &(m\circ g)\circ h\ar[r]^\alpha & m\circ (g\circ h)\ar[r]^{m*\gamma}&m\circ(k\circ f)\ar[r]^{\alpha^{-1}}& (m\circ k)\circ f})
$$
as in the following diagram
$$
\xymatrix{
A\ar[r]^h\ar[d]_f\ar@{}[dr]|\gamma & C\ar[d]_{g}\ar[r]^l\ar@{}[dr]|\delta& E\ar[d]^i
\\
B\ar[r]_k& D\ar[r]_m& F.}
$$
  Given 2-cells
  $(\varphi_{0},\varphi_{1})\colon(h,k,\gamma)\to(h',k',\gamma')$ and
  $(\psi_{0},\psi_{1})\colon(l,m,\delta)\to(l',m',\delta')$ we define
  the horizontal composite by
  \begin{align*}
    (\psi_{0},\psi_{1})*(\varphi_{0},\varphi_{1}) & := (\psi_{0}*\varphi_{0},\psi_{1}*\varphi_{1}).
  \end{align*}
\end{description}
The verification that, with these definitions, $\CC^{\to}$ forms a
bicategory is lengthy, but straightforward, and is left to the reader.  Note that when $\CC$ is a 2-category, so is $\CC^{\to}$.

\subsection{The 2-categorical case}\label{section:pseudofunctors}

We later will be concerned with pseudofunctors $F\colon \CC^{\myop}\to\cat$ where
$\CC$ is a category understood as having a trivial 2-category
structure.  This means precisely that for each object $U$ of $\CC$
there is a category $F(U)$ and for each map $f\colon V\to U$ there is
a functor $F(f)\colon F(U)\to F(V)$.  We will often denote the action of
$F(f)$ on $x\in F(U)$ by $x\cdot f$ or, when the map $f$ is
understood, by $x|_{V}$.  For each object $U$ we require a
distinguished natural isomorphism $\upsilon_{U}\colon F(1_{U})\to 1_{F(U)}$.  Finally, for
$f\colon V\to U$ and $g\colon  E\to V$, we require a distinguished natural
isomorphism $\varphi_{f,g}\colon F(g\circ f)\to F(g)\circ F(f)$.
When $f$ and $g$ are understood we omit subscripts and simply write
$\varphi$.  Similarly, we sometimes write $\upsilon$ instead of $\upsilon_{U}$.

Assume given a fixed object $U$ of $\CC$ together with $x$ in $F(U)$
and arrows $f\colon U_{\alpha}\to U$ and $g\colon U_{\beta}\to U$.  In this
situation we often denote by $x|_{\alpha}$ the object $x|_{U_{\alpha}}$ and,
similarly, we denote by $x|_{\alpha\beta}$ the object
$x|_{U_{\alpha}\times_{U}U_{\beta}}$.  In this situation, we will make
use below of the isomorphism from $x|_{\alpha}|_{\alpha\beta}$ to
$x|_{\beta}|_{\alpha\beta}$ constructed using the coherence maps $\varphi$ and
for which we introduce the notation $\sigma_{\beta\alpha}(x)$.
Explicitly, $\sigma_{\beta\alpha}(x)$ is defined to be the composite
\begin{align*}
  \xy
  {\ar^{\varphi^{-1}(x)}(0,0)*+{x|_{\alpha}|_{\alpha\beta}};(30,0)*+{x|_{\alpha\beta}}};
  {\ar^{\varphi(x)}(30,0)*+{x|_{\alpha\beta}};(60,0)*+{x|_{\beta}|_{\alpha\beta}.}};
  \endxy
\end{align*}
We also remark that $\sigma_{\alpha\beta}$ is the inverse of
$\sigma_{\beta\alpha}$.  We similarly write
$\sigma_{\alpha,\beta\gamma}(x)$ for the map
$x|_{\alpha\gamma}|_{\alpha\beta\gamma}\to
x|_{\alpha\beta}|_{\alpha\beta\gamma}$ which is defined in the same
way as the composite
\begin{align*}
   \xy
  {\ar^{\varphi^{-1}(x)}(0,0)*+{x|_{\alpha\gamma}|_{\alpha\beta\gamma}};(30,0)*+{x|_{\alpha\beta\gamma}}};
  {\ar^{\varphi(x)}(30,0)*+{x|_{\alpha\beta\gamma}};(60,0)*+{x|_{\alpha\beta}|_{\alpha\beta\gamma},}};
  \endxy
\end{align*}
for $x$ an object of $F(U_{\alpha})$.

\section{Fibrant objects and fibration structures}\label{section:fibration_structures}

We will now axiomatize two bicategorical notions: bicategories with
systems of fibrant objects and bicategories with fibration
structures.  The former suffices for the construction of the homotopy
bicategory.  However, the latter concept, which is a refinement of
the former, captures additional structure present in certain examples
and provides additional
structure such as path objects for the homotopy category.

\subsection{Systems of fibrant objects}

We will now turn to consider the axiomatic structure on a bicategory
which will allow us to form the homotopy bicategory and prove that it
possesses the correct universal property.
\begin{definition}
  For arrows $f\colon A\to B$ and $g\colon C\to D$ in a bicategory $\CC$, we write $f\pitchfork g$ to
  indicate that for any square of the form
  \begin{align}\label{eq:lift_diagram}
    \xy
    {\ar^{h}(0,20)*+{A};(20,20)*+{C}};
    {\ar_{f}(0,20)*+{A};(0,0)*+{B}};
    {\ar^{g}(20,20)*+{C};(20,0)*+{D}};
    {\ar_{k}(0,0)*+{B};(20,0)*+{D}};
    {\ar^{\gamma}@{=>}(10,13)*+{};(10,7)*+{}};
    \endxy
  \end{align}
  with $\gamma$ an invertible 2-cell, there exists a map $l\colon B\to
  C$ together with invertible 2-cells $\lambda\colon h\iso l\circ f$ and
  $\rho\colon g\circ l\iso k$ such that
  \begin{align*}
    \xy
    {\ar_{\gamma}(0,20)*+{g\circ h};(0,0)*+{k\circ f}};
    {\ar^(.45){g*\lambda}(0,20)*+{g\circ h};(20,20)*+{g\circ(l\circ f)}};
    {\ar^{\alpha^{-1}}(20,20)*+{g\circ (l\circ f)};(20,0)*+{(g\circ
        l)\circ f}};
    {\ar^{\rho*f}(20,0)*+{(g\circ l)\circ f};(0,0)*+{k\circ f}};
    \endxy
  \end{align*}
  commutes in $\CC(A,D)$.

  Given a class $\frak{M}$ of maps in $\CC$ we write
  $\frak{M}\pitchfork g$ to indicate that $f\pitchfork g$ for all $f$
  in $\frak{M}$.  For $C$ an object of $\CC$, we write
  $\frak{M}\pitchfork C$ to indicate that, for all maps $f\colon A\to
  B$ and $h\colon A\to C$, if $f$ is in $\frak{M}$, then there exists
  a map $l\colon B\to C$ and an invertible 2-cell $h\Rightarrow l\circ f$.
\end{definition}
Observe that if a bicategory $\CC$ has a terminal object $1$, then
$\frak{M}\pitchfork C$ if and only if \mbox{$\frak{M}\pitchfork ( C\to 1)$}.
\begin{definition}\label{def:sfo}
  A \myemph{system of fibrant objects} on a bicategory $\CC$ consists
  of a collection of maps $\frak{W}$ (\emph{weak equivalences}) of
  $\CC$ together with a pseudofunctor $Q\colon\CC\to\CC$
  (\emph{fibrant replacement}) and a pseudonatural transformation
  $\eta\colon 1_{\CC}\to Q$ such that the following axioms are satisfied:
  \begin{description}
  \item[Identities] All identity arrows $1_{A}\colon A\to A$ are in $\frak{W}$.
  \item[3-for-2] Given a diagram
    \begin{align*}
      \xy
      {\ar^{f}(0,0)*+{A};(24,0)*+{C}};
      {\ar_{g}(0,0)*+{A};(12,-14)*+{B}};
      {\ar_{h}(12,-14)*+{B};(24,0)*+{C}};
      {\ar@{=>}^{\gamma}(12,-3)*+{};(12,-8)*+{}};
      \endxy
    \end{align*}
    with $\gamma$ an isomorphism, if any two of $f,g$ and $h$ are weak
    equivalences, then so is the third.
  \item[Fibrant Replacement] The components
    of $\eta$ are weak equivalences and $\frak{W}\pitchfork Q(A)$ for
    any object $A$ of $\CC$.
  \end{description}
\end{definition}
The notion of a fibration structure on a bicategory $\CC$ is a slight
refinement of the notion of a system of fibrant objects:
\begin{definition}\label{def:fs}
  A \myemph{fibration structure} on a bicategory $\CC$ with terminal
  object $1$ is given by collections of maps $\frak{W}$ (\emph{weak
    equivalences}) and $\frak{F}$
  (\emph{fibrations}) of $\CC$ such that $\frak{W}$ satisfies the
  identities and 3-for-2 conditions from Definition \ref{def:sfo}
  above and such that the following additional
  axioms are satisfied:
  \begin{description}
  \item[Lifting] $p\colon E\to B$ is a fibration if and only if $\frak{W}\pitchfork p$.
  \item[Factorization] There exists a pseudofunctor
    $Q\colon\CC^{\to}\to\CC^{\to}$ together with a pseudonatural
    transformation $\eta\colon 1_{\CC^{\to}}\to Q$ such that
    $\partial_{1}\circ Q=\partial_{1}$,
    $\partial_{1}\circ\eta=\partial_{1}$, and, for each
    $f\colon A\to B$ in $\CC$, the arrow part of $\eta_{f}$ is a weak
    equivalence and $Q(f)$ is a fibration.  Here
    $\partial_{1}$ is the pseudofunctor $\CC^{\to}\to\CC$ which
    projects onto the codomain.
  \end{description}
\end{definition}
Note that every fibration structure on a bicategory $\CC$ determines a
corresponding system of fibrant objects.

\begin{remark}
 When we apply the factorization condition to the diagonal $\Delta_A\colon A\rightarrow A\times A$,
we obtain a diagram
$$
\xymatrix{
A\ar[d]_{\Delta_A}\ar@{}[dr]|\cong\ar[r]^{\eta_{\Delta_A}}& \partial_0 Q\Delta_A\ar[d]^{Q\Delta_A}
\\
A\times A\ar@{=}[r] & A\times A
}$$
Here, $\eta_{\Delta_A}$ is a weak equivalence and $Q\Delta_A$ is a fibration. So we find that we can take $A^I=\partial_0 Q\Delta_A$
as a path object for $A$ and the classical factorization lemma holds up to an invertible 2-cell.
Furthermore, when we take $d_i=\pi_i\circ Q\Delta_A\colon A^I\rightarrow A$ we obtain a fibration with the property that 
$d_i\circ \eta_{\Delta_A}\cong1_A$.
\end{remark}

For the remainder of this section we assume that we are working in a bicategory $\CC$ with
a system of fibrant objects.
\begin{definition}
  An object $A$ of $\CC$ is \myemph{fibrant} when $\frak{W}\pitchfork A$.
\end{definition}
\begin{lemma}
  \label{lemma:invertible_general}
  If $f\colon A\to B$ is a weak equivalence between fibrant objects, then $f$ is an
  equivalence.
  \begin{proof}
    First, since $A$ is fibrant there exists a map $f'\colon B\to A$ and an
    invertible 2-cell 
    \begin{align*}
      \xy
      {\ar^{1_{A}}(0,15)*+{A};(20,15)*+{A.}};
      {\ar_{f}(0,15)*+{A};(0,0)*+{B}};
      {\ar_{f'}(0,0)*+{B};(20,15)*+{A.}};
      {\ar@{=>}(7,13)*+{};(7,7)*+{}};
      \endxy
    \end{align*}
    It follows from the 3-for-2 property that $f'$ is also a weak
    equivalence.  Therefore, since $B$ is fibrant, there exists
    another map $f''\colon A\to B$ and an invertible 2-cell 
    \begin{align*}
      \xy
      {\ar^{1_{B}}(0,15)*+{B};(20,15)*+{B.}};
      {\ar_{f'}(0,15)*+{B};(0,0)*+{A}};
      {\ar_{f''}(0,0)*+{A};(20,15)*+{B.}};
      {\ar@{=>}(7,13)*+{};(7,7)*+{}};
      \endxy
    \end{align*}
    Now, the 2-cells above, together with the coherence 2-cells of
    $\CC$, give us an isomorphism $f\iso f''$ and therefore
    $f'$ is the pseudo-inverse of $f$, as required.
  \end{proof}
\end{lemma}
\begin{definition}
  The \myemph{homotopy bicategory $\hocat(\CC)$ of $\CC$} is the
  full sub-bicategory of fibrant objects of $\CC$.
\end{definition}
We denote by $I\colon \CC\to\hocat(\CC)$ the pseudofunctor induced by
$Q\colon \CC\to\CC$.  It is an immediate consequence of Lemma \ref{lemma:invertible_general}
that $I$ sends weak equivalences to equivalences.  For any bicategory
$\DD$, let $[\CC,\DD]_{\frak{W}}$ denote the
sub-bicategory of $[\CC,\DD]$ consisting of those pseudofunctors which
send maps in $\frak{W}$ to equivalences.  Let
$J\colon\hocat(\CC)\to\CC$ be the inclusion and observe that $Q=J\circ I$.

We will now prove that $I$ is the universal map from $\CC$
to a bicategory which sends weak equivalences to equivalences.
\begin{theorem}\label{theorem:main}
  For any bicategory $\DD$, $I$ induces an equivalence of bicategories
  \begin{align*}
    \xy
    {\ar^{[I,\DD]}(0,0)*+{[\hocat(\CC),\DD]};(35,0)*+{[\CC,\DD]_{\frak{W}},}};
    \endxy
  \end{align*}
  where the subscript
  $\frak{W}$ indicates that we are considering only those pseudofunctors
  which send weak equivalences to equivalences.
  \begin{proof}
    Precomposition with $J$ gives a pseudofunctor
    $[\CC,\DD]_{\frak{W}}\to [\hocat(\CC),\DD]$ which we denote by $[J,\DD]$.  
    The pseudonatural transformation $\eta\colon 1_{\hocat(\CC)}\to
    I\circ J$ (obtained by restricting $\eta$ to $\hocat(\CC)$) induces a pseudonatural transformation
    $[\eta,\DD]\colon 1_{[\hocat(\CC),\DD]}\to[J,\DD]\circ[I,\DD]$.
    Observe that, by Lemmas \ref{lemma:pseudonat} and
    \ref{lemma:invertible_general}, $\eta\colon 1_{\hocat(\CC)}\to
    I\circ J$ is an equivalence.  Therefore the induced $[\eta,\DD]$ is also an
    equivalence.

    On the other hand, for $F$ in $[\CC,\DD]_{\frak{W}}$, Lemma
    \ref{lemma:pseudonat} exhibits $F\eta\colon F\to FQ$ as an adjoint
    equivalence.  Let $\vartheta^{F}$ denote the adjoint pseudoinverse
    of $F\eta$.  Allowing $F$ to vary, we have that $\vartheta$ is an
    equivalence $[I,\DD]\circ[J,\DD]\to 1_{[\CC,\DD]_{\frak{W}}}$.
  \end{proof}
\end{theorem}

\section{Stacks and local fibrations}\label{section:local_fib}

We will now begin developing the machinery required to explain our
first example of a fibration system in a bicategory (actually, in this
case a 2-category): the 2-category of prestacks.  In this section we
recall some of the basic notions involved and we also introduce a
fibered version of the usual category of descent data that will allow
us to describe the maps, which we call \emph{local fibrations}, that
provide the fibrations in the fibration structure for the 2-category of prestacks.

\subsection{Coverings and sites}

Throughout we assume given a fixed site $(\CC,J)$ for $\CC$ a category
with finite limits.  Given an object
$U$ of $\CC$, recall that a \emph{sieve on $U$} is a family of maps
with codomain $U$ which is a right ideal for composition.  To say that
$(\CC,J)$ is a site then means that $J$  assigns to each object $U$ of
$\CC$ a collection $J(U)$ of sieves on $U$ (called
\emph{covering sieves}, \emph{covering families} or \emph{covers}) in such a way that the following conditions
are satisfied:
\begin{enumerate}
\item The maximal sieve on $U$, which consists of all arrows with
  codomain $U$, is in $J(U)$.
\item For a cover $\SS$ in $J(U)$ and $g\colon V \to U$, the
  sieve \mbox{$g^{*}(\SS):=\{ f \colon E\to V \;|\; g\circ f \in \SS\}$} is in $J(V)$.
\item Given $\SS$ in $J(U)$ and a sieve $\RR$ on $U$, if
  $f^{*}(\RR)$ is in $J(V)$ for all \mbox{$f\colon V\to U$} in $\SS$, then $\RR$
  is also in $J(U)$.
\end{enumerate}
We will sometimes also work with the notion of a \emph{basis} for
covers.  A basis consists of an operation $K$ which assigns to objects
$U$ of $\CC$ a collection $K(U)$ of families of maps with codomain $U$ such that
\begin{enumerate}
\item The singleton family $(f\colon U'\to U)$ is in $K(U)$ when $f$
  is an isomorphism.
\item If $U_{\alpha}\to U$ is in $K(U)$ and $V\to U$ is any map, then
  $V\times_{U}U_{\alpha}\to V$ is in $K(V)$.
\item If $(f_{\alpha}\colon U_{\alpha}\to U)$ is a family of maps in
  $K(U)$ and $(f^{\alpha}_{\beta}\colon U^{\alpha}_{\beta}\to
  U_{\alpha})$ is in $K(U_{\alpha})$, then the family of maps $(f_{\alpha}\circ
  f^{\alpha}_{\beta})$ is in $K(U)$.
\end{enumerate}
If $K$ is a basis, then $K$ generates a site $(\CC,J)$ by
letting $S$ be in $J(U)$ if and only if there exists a family $R$ in
$K(U)$ such that $R\subseteq S$.

Readers unfamiliar with sites and sheaves may consult
\cite{MacLane:SGL}.

\subsection{Stacks}

Given a site $(\CC,J)$, a pseudofunctor
$F\colon \CC^{\myop}\to\cat$ is a \myemph{stack} when, for any
cover $(f_{\alpha}\colon U_{\alpha}\to U)_{\alpha}$ of $U$, the canonical map
\begin{align}\label{eq:is_stack}
  F(U)\to \limit_{\alpha}F(U_{\alpha})
\end{align}
is a weak equivalence of categories.  Note that here
$\limit_{\alpha}F(U_{\alpha})$ indicates the
\emph{pseudolimit} and \emph{not} the strict limit (an elementary
description can be found below for which we refer the reader to
Example \ref{example:stacks}).  $\stacks(\CC)$
denotes the full subcategory of $\pseudo{\CC^{\myop}}{\cat}$
consisting of stacks.  

Note that if a basis $K$ generates the covering sieves of a site
$(\CC,J)$, then it suffices, in order to tell whether $F$ is a stack,
to test on the families of maps $\UU$ in $K(U)$.

\subsection{Prestacks}\label{section:prestacks}

Given a pseudofunctor $F\colon \CC^{\myop}\to\cat$ and objects $a$ and
$b$ of $F(U)$, there is an induced functor
\begin{align*}
  \xy
  {\ar^(.6){F(a,b)}(0,0)*+{(\CC/U)^{\myop}};(30,0)*+{\sets}};
  \endxy
\end{align*}
which sends
\begin{align*}
  \xy
  {\ar@{|->}(0,0)*+{f\colon V\to U};(30,0)*+{\hom_{FV}(a|_{V},b|_{V}).}};
  \endxy
\end{align*}
A pseudofunctor $F\colon \CC^{\myop}\to\cat$ is a \myemph{prestack} if,
for any $U$ and $a,b\in F(U)$, $F(a,b)$ is a sheaf.  In particular, in
a prestack it is possible to construct arrows in the categories $F(U)$
locally (i.e., on a cover). Note that every stack is a prestack.

\subsection{Descent data}\label{section:descent_data}

The following definition generalizes the familiar definition of the
category of descent data by making this data vary relative to a fixed
morphism.
\begin{definition}
  Given $p\colon E\to B$ in $\pseudo{\CC^{\myop}}{\cat}$ and a
  cover $\SS=(f_{\alpha}\colon U_{\alpha}\to U)_{\alpha}$ of some $U$ we define
  the category $\descent{p}{\SS}$ as follows
  \begin{description}
  \item[Objects] An object is a tuple
    $(b,(e_{\alpha}),(\psi_{\alpha}),(\vartheta_{\alpha\beta}))$ where
    $b\in B(U)$, $e_{\alpha}\in E(U_{\alpha})$, $\psi_{\alpha}$ is an
    isomorphism $p(e_{\alpha})\to b|_{\alpha}$ and
    $\vartheta_{\alpha\beta}$ is an isomorphism
    $e_{\beta}|_{\alpha\beta}\to e_{\alpha}|_{\alpha\beta}$.  This
    data is furthermore required to satisfy the conditions that 
    \begin{align*}
      \vartheta_{\alpha\alpha} & = 1_{e_{\alpha}}
    \end{align*}
    and that the diagrams
    \begin{align*}
      \xy
      {\ar_{\sigma_{\gamma,\alpha\beta}(e_{\gamma})}(30,20)*+{e_{\gamma}|_{\beta\gamma}|_{\alpha\beta\gamma}};(0,20)*+{e_{\gamma}|_{\alpha\gamma}|_{\alpha\beta\gamma}}};
      {\ar_{\vartheta_{\alpha\gamma}|_{\alpha\beta\gamma}}(0,20)*+{e_{\gamma}|_{\alpha\gamma}|_{\alpha\beta\gamma}};(0,0)*+{e_{\alpha}|_{\alpha\gamma}|_{\alpha\beta\gamma}}};
      {\ar_{\sigma_{\alpha,\beta\gamma}(e_{\alpha})}(0,0)*+{e_{\alpha}|_{\alpha\gamma}|_{\alpha\beta\gamma}};(30,0)*+{e_{\alpha}|_{\alpha\beta}|_{\alpha\beta\gamma}}};
      {\ar^{\vartheta_{\beta\gamma}|_{\alpha\beta\gamma}}(30,20)*+{e_{\gamma}|_{\beta\gamma}|_{\alpha\beta\gamma}};(60,20)*+{e_{\beta}|_{\beta\gamma}|_{\alpha\beta\gamma}}};
      {\ar^{\sigma_{\beta,\alpha\gamma}(e_{\beta})}(60,20)*+{e_{\beta}|_{\beta\gamma}|_{\alpha\beta\gamma}};(60,0)*+{e_{\beta}|_{\alpha\beta}|_{\alpha\beta\gamma}}};
      {\ar^{\vartheta_{\alpha\beta}|_{\alpha\beta\gamma}}(60,0)*+{e_{\beta}|_{\alpha\beta}|_{\alpha\beta\gamma}};(30,0)*+{e_{\alpha}|_{\alpha\beta}|_{\alpha\beta\gamma}}};
      \endxy
    \end{align*}
    and
    \begin{align*}
      \xy
      {\ar^{\eta(e_{\beta})}(0,20)*+{p(e_{\beta}|_{\alpha\beta})};(25,20)*+{p(e_{\beta})|_{\alpha\beta}}};
      {\ar^{\psi_{\beta}|_{\alpha\beta}}(25,20)*+{p(e_{\beta})|_{\alpha\beta}};(50,20)*+{b|_{\beta}|_{\alpha\beta}}};
      {\ar_{p(\vartheta_{\alpha\beta})}(0,20)*+{p(e_{\beta}|_{\alpha\beta})};(0,0)*+{p(e_{\alpha}|_{\alpha\beta})}};
      {\ar_{\eta(e_{\alpha})}(0,0)*+{p(e_{\alpha}|_{\alpha\beta})};(25,0)*+{p(e_{\alpha})|_{\alpha\beta}}};
      {\ar_{\psi_{\alpha}|_{\alpha\beta}}(25,0)*+{p(e_{\alpha})|_{\alpha\beta}};(50,0)*+{b|_{\alpha}|_{\alpha\beta}}};
      {\ar^{\sigma_{\alpha\beta}(b)}(50,20)*+{b|_{\beta}|_{\alpha\beta}};(50,0)*+{b|_{\alpha}|_{\alpha\beta}}};
      \endxy
    \end{align*}
    commute, where the $\eta$ are the natural isomorphisms associated
    to $p$.
  \item[Arrows] An arrow
    $(b,(e_{\alpha}),(\psi_{\alpha}),(\vartheta_{\alpha\beta}))\to(b',(e'_{\alpha}),(\psi'_{\alpha}),(\vartheta'_{\alpha\beta}))$
    is given by a pair $(g,(g_{\alpha}))$ such that $g\colon b\to b'$ in
    $B(U)$ and $g_{\alpha}\colon e_{\alpha}\to e'_{\alpha}$ in
    $E(U_{\alpha})$. This data is subject to the requirements that the
    diagrams
    \begin{align*}
      \xy
      {\ar^{\vartheta_{\alpha\beta}}(0,20)*+{e_{\beta}|_{\alpha\beta}};(20,20)*+{e_{\alpha}|_{\alpha\beta}}};
      {\ar_{g_{\beta}|_{\alpha\beta}}(0,20)*+{e_{\beta}|_{\alpha\beta}};(0,0)*+{e'_{\beta}|_{\alpha\beta}}};
      {\ar^{g_{\alpha}|_{\alpha\beta}}(20,20)*+{e_{\alpha}|_{\alpha\beta}};(20,0)*+{e'_{\alpha}|_{\alpha\beta}}};
      {\ar_{\vartheta'_{\alpha\beta}}(0,0)*+{e'_{\beta}|_{\alpha\beta}};(20,0)*+{e'_{\alpha}|_{\alpha\beta}}};
      \endxy
    \end{align*}
    and 
    \begin{align*}
      \xy
      {\ar^{\psi_{\alpha}}(0,20)*+{p(e_{\alpha})};(20,20)*+{b|_{\alpha}}};
      {\ar^{g|_{\alpha}}(20,20)*+{b|_{\alpha}};(20,0)*+{b'|_{\alpha}}};
      {\ar_{p(g_{\alpha})}(0,20)*+{p(e_{\alpha})};(0,0)*+{p(e'_{\alpha})}};
      {\ar_{\psi'_{\alpha}}(0,0)*+{p(e'_{\alpha})};(20,0)*+{b'|_{\alpha}}};
      \endxy
    \end{align*}
    commute.
  \end{description}
\end{definition}
There is a projection functor \mbox{$\pi\colon \descent{p}{\SS}\to B(U)$}.
When $p$ is the canonical map $E\to 1$ into the terminal object we
write $\descent{E}{\SS}$ instead of $\descent{p}{\SS}$ and observe
that this is the usual category of descent data.  There is also an evident functor
\begin{align*}
  \xy
  {\ar^{\Phi_{\SS}}(0,15)*+{E(U)};(30,15)*+{\descent{p}{\SS}}};
  {\ar_{p}(0,15)*+{E(U)};(15,0)*+{B(U)}};
  {\ar^{\pi}(30,15)*+{\descent{p}{\SS}};(15,0)*+{B(U)}};
  \endxy
\end{align*}
which sends an object $e$ of $E(U)$ to the tuple
\begin{align*}
  \bigl(p(e),(e|_{\alpha}),(1_{p(e|_{\alpha})}),(\sigma_{\alpha\beta}(e))\bigr)
\end{align*}
where $\sigma_{\alpha\beta}$ is as in Section \ref{section:pseudofunctors}.

Given a commutative triangle
\begin{align*}
  \xy
  {\ar^{f}(0,15)*+{E};(30,15)*+{E'}};
  {\ar_{p}(0,15)*+{E};(15,0)*+{B}};
  {\ar^{p'}(30,15)*+{E'};(15,0)*+{B}};
  \endxy
\end{align*}
in $\pseudo{\CC^{\myop}}{\cat}$ and a covering family $\SS$ of $U$,
there is a corresponding commutative diagram
\begin{align*}
  \xy
  {\ar^{f_{*}}(0,15)*+{\descent{p}{\SS}};(30,15)*+{\descent{p'}{\SS}}};
  {\ar_{\pi}(0,15)*+{\descent{p}{\SS}};(15,0)*+{B(U)}};
  {\ar^{\pi}(30,15)*+{\descent{p'}{\SS}};(15,0)*+{B(U)}};
  \endxy
\end{align*}
in $\cat$.  Here the functor $f_{*}$ acts as follows:
\begin{description}
\item[On objects] $f_{*}$ sends 
  $(b,(e_{\alpha}),(\psi_{\alpha}),(\vartheta_{\alpha\beta}))$ in
  $\descent{p}{\UU}$ to
  $(b,(f(e_{\alpha})),(\psi_{\alpha}),\xi_{\alpha\beta})$ where
  $\xi_{\alpha\beta}$ is the composite
  \begin{align*}
    \xy
    {\ar(0,0)*+{f(e_{\beta})|_{\alpha\beta}};(25,0)*+{f(e_{\beta}|_{\alpha\beta})}};
    {\ar^{f(\vartheta_{\alpha\beta})}(25,0)*+{f(e_{\beta}|_{\alpha\beta})};(50,0)*+{f(e_{\alpha}|_{\alpha\beta})}};
    {\ar(50,0)*+{f(e_{\alpha}|_{\alpha\beta})};(75,0)*+{f(e_{\alpha})|_{\alpha\beta}}};
    \endxy
  \end{align*}
  where the unnamed arrows are from the coherence isomorphisms
  associated to the pseudonatural transformation $f$.
\item[On arrows] An arrow $(g,(g_{\alpha}))$ in
  $\descent{p}{\UU}$ is sent to $(g,(f(g_{\alpha})))$.
\end{description}
Moreover, this construction is functorial in the sense that
$\descent{-}{\SS}$ is a functor from $\pseudo{\CC^{\myop}}{\cat}/B$
to $\cat/B(U)$.  This fact is a special case of a more general result
to which we now turn.
\begin{lemma}
  \label{lemma:functorial_squares}
  For a fixed object $U$ of $\CC$ and a covering family $\SS$ of $U$,
  $\descent{-}{\SS}$ is a functor $(\pseudo{\CC^{\myop}}{\cat})^{\rightarrow}\to\cat$.
  \begin{proof}
    Given a square 
    \begin{align}\label{eq:iso_square}
      \xy
      {\ar^{h}(0,20)*+{A};(20,20)*+{E}};
      {\ar_{i}(0,20)*+{A};(0,0)*+{C}};
      {\ar^{p}(20,20)*+{E};(20,0)*+{B}};
      {\ar_{k}(0,0)*+{C};(20,0)*+{B}};
      {\ar^{\gamma}@{=>}(10,13)*+{};(10,7)*+{}};
      \endxy
    \end{align}
    in $\pseudo{\CC^{\myop}}{\cat}$ with $\gamma$ invertible, the induced functor
    $(h,k,\gamma)_{*}\colon \descent{i}{\SS}\to\descent{p}{\SS}$ sends descent data
    $(c,(a_{\alpha}),(\psi_{\alpha}),(\vartheta_{\alpha\beta}))$ in
    $\descent{i}{\UU}$ to 
    $(k(c),(h(a_{\alpha})),(\hat{\psi}_{\alpha}),(\hat{\vartheta}_{\alpha\beta}))$
    in $\descent{p}{\UU}$
    where $\hat{\psi}_{\alpha}$ is the composite
    \begin{align*}
      \xy
      {\ar^{\gamma(a_{\alpha})}(0,0)*+{p(h(a_{\alpha}))};(25,0)*+{k(i(a_{\alpha}))}};
      {\ar^{k(\psi_{\alpha})}(25,0)*+{k(i(a_{\alpha}))};(50,0)*+{k(c|_{\alpha})}};
      {\ar(50,0)*+{k(c|_{\alpha})};(75,0)*+{k(c)|_{\alpha}}};
      \endxy
    \end{align*}
    and $\hat{\vartheta}_{\alpha\beta}$ is the composite
    \begin{align*}
      \xy
      {\ar(0,0)*+{h(a_{\beta})|_{\alpha\beta}};(25,0)*+{h(a_{\beta}|_{\alpha\beta})}};
      {\ar^{h(\vartheta_{\alpha\beta})}(25,0)*+{h(a_{\beta}|_{\alpha\beta})};(50,0)*+{h(a_{\alpha}|_{\alpha\beta})}};
      {\ar(50,0)*+{h(a_{\alpha}|_{\alpha\beta})};(75,0)*+{h(a_{\alpha})|_{\alpha\beta}.}};
      \endxy
    \end{align*}
    Here the unnamed arrows are from the coherence isomorphisms
    associated to the pseudonatural transformations.
  \end{proof}
\end{lemma}
Observe that, given a square (\ref{eq:iso_square}) in
$\pseudo{\CC^{\myop}}{\cat}$ and a cover $\SS$ of some $U$,
the following diagram commutes:
\begin{align*}
  \xy
  {\ar^{(h,k,\gamma)_{*}}(0,15)*+{\descent{i}{\SS}};(30,15)*+{\descent{p}{\SS}}};
  {\ar_{\pi}(0,15)*+{\descent{i}{\SS}};(0,0)*+{C(U)}};
  {\ar^{\pi}(30,15)*+{\descent{p}{\SS}};(30,0)*+{B(U).}};
  {\ar_{k}(0,0)*+{C(U)};(30,0)*+{B(U).}};
  \endxy
\end{align*}
On the other hand, we merely have a natural isomorphism
$\hat{\gamma}$ as indicated in the following diagram:
\begin{align*}
  \xy
  {\ar^{h}(0,15)*+{A(U)};(30,15)*+{E(U)}};
  {\ar_{\Phi_{\SS}}(0,15)*+{A(U)};(0,0)*+{\descent{i}{\SS}}};
  {\ar^{\Phi_{\SS}}(30,15)*+{E(U)};(30,0)*+{\descent{p}{\SS}}};
  {\ar_{(h,k,\gamma)_{*}}(0,0)*+{\descent{i}{\SS}};(30,0)*+{\descent{p}{\SS}}};
  {\ar@{=>}^{\hat{\gamma}}(15,10)*+{};(15,5)*+{}};
  \endxy
\end{align*}
which, for $a$ an object of $A(U)$, is the map of descent data
\begin{align*}
  (\gamma(a)\colon ph(a)\to ki(a),(h(a)|_{\alpha}\to
  h(a|_{\alpha})))\colon \Phi_{\SS}(h(a))\to (h,k,\gamma)_{*}(\Phi_{\SS}(a)).
\end{align*}
This has the property that
\begin{align*}
  \xy
  {\ar^{h}(0,15)*+{A(U)};(30,15)*+{E(U)}};
  {\ar_{\Phi_{\SS}}(0,15)*+{A(U)};(0,0)*+{\descent{i}{\SS}}};
  {\ar^{\Phi_{\SS}}(30,15)*+{E(U)};(30,0)*+{\descent{p}{\SS}}};
  {\ar_{(h,k,\gamma)_{*}}(0,0)*+{\descent{i}{\SS}};(30,0)*+{\descent{p}{\SS}}};
  {\ar@{=>}^{\hat{\gamma}}(15,10)*+{};(15,5)*+{}};
  {\ar(0,0)*+{\descent{i}{\SS}};(30,0)*+{\descent{p}{\SS}}};
  {\ar_{\pi}(0,0)*+{\descent{i}{\SS}};(0,-15)*+{C(U)}};
  {\ar^{\pi}(30,0)*+{\descent{p}{\SS}};(30,-15)*+{B(U)}};
  {\ar_{k}(0,-15)*+{C(U)};(30,-15)*+{B(U)}};
  {(45,0)*+{=}};
  {\ar^{h}(60,7.5)*+{A(U)};(90,7.5)*+{E(U)}};
  {\ar_{i}(60,7.5)*+{A(U)};(60,-7.5)*+{C(U)}};
  {\ar^{p}(90,7.5)*+{E(U)};(90,-7.5)*+{B(U)}};
  {\ar_{k}(60,-7.5)*+{C(U)};(90,-7.5)*+{B(U)}};
  {\ar@{=>}^{\gamma}(75,3)*+{};(75,-2)*+{}};
  \endxy
\end{align*}

The construction of the category of descent data is also functorial in
the second argument in the sense that if $\SS$ and $\RR$ are both
covers of some $U$ with $\RR\subseteq\SS$, then there exists an
associated restriction functor
$\cdot|_{\RR}\colon \descent{p}{\SS}\to\descent{p}{\RR}$ which acts by
restricting descent data to those maps in $\RR$.  These restrictions
satisfy the functoriality condition
$(\cdot|_{\UU})\circ(\cdot|_{\RR})=\cdot|_{\UU}$ and are
well-behaved with respect to the associated maps
$\Phi_{\SS}\colon E(U)\to\descent{p}{\SS}$, in the sense that the diagram
\begin{align*}
  \xy
  {\ar^{\Phi_{\SS}}(0,15)*+{E(U)};(30,15)*+{\descent{p}{\SS}}};
  {\ar_{\Phi_{\RR}}(0,15)*+{E(U)};(15,0)*+{\descent{p}{\RR}}};
  {\ar^{\cdot|_{\RR}}(30,15)*+{\descent{p}{\SS}};(15,0)*+{\descent{p}{\RR}}};
  \endxy
\end{align*}
commutes for any $\RR\subseteq\SS$.

In addition to the functorial behavior of $\descent{-}{-}$ described
above, if we are given a fixed $p\colon E\to B$, a cover $\UU$ of $U$ and a
map $g\colon V\to U$ in the site, we obtain a further restriction functor
$g^{*}\colon \descent{p}{\UU}\to\descent{p}{g^{*}(\UU)}$ which sends descent
data $(b,(e_{\alpha}),(\psi_{\alpha}),(\vartheta_{\alpha\beta}))$ to
the descent data given by:
\begin{itemize}
\item the object $b|_{V}$ of $B(V)$;
\item the family of objects $(e_{\alpha})$ (this makes sense by virtue
  of the definition of $g^{*}(\UU)$);
\item the family of maps given by the composites
  \begin{align*}
    \xy
    {\ar^{\psi_{\alpha}}(0,0)*+{p(e_{\alpha})};(20,0)*+{b|_{\alpha}}};
    {\ar(20,0)*+{b|_{\alpha}};(40,0)*+{(b|_{V})|_{\alpha},}};
    \endxy 
  \end{align*}
  which we denote by $g^{*}(\psi)_{\alpha}$ when no confusion will
  result; and
\item the family of maps given by the composites 
  \begin{align*}
    \xy
    {\ar(0,0)*+{e_{\beta}|_{U_{\alpha}\times_{V}U_{\beta}}};(30,0)*+{e_{\beta}|_{U_{\alpha}\times_{U}U_{\beta}}|_{U_{\alpha}\times_{V}U_{\beta}}}};
    {\ar^{\vartheta_{\alpha\beta}|_{U_{\alpha}\times_{V}U_{\beta}}}(30,0)*+{e_{\beta}|_{U_{\alpha}\times_{U}U_{\beta}}|_{U_{\alpha}\times_{V}U_{\beta}}};(70,0)*+{e_{\alpha}|_{U_{\alpha}\times_{U}U_{\beta}}|_{U_{\alpha}\times_{V}U_{\beta}}}};
    {\ar(70,0)*+{e_{\alpha}|_{U_{\alpha}\times_{U}U_{\beta}}|_{U_{\alpha}\times_{V}U_{\beta}}};(100,0)*+{e_{\alpha}|_{U_{\alpha}\times_{V}U_{\beta}}}};
    \endxy
  \end{align*}
  where $\vartheta_{\alpha\beta}$ is here restricted along the induced map $U_{\alpha}\times_{V}U_{\beta}\to
  U_{\alpha}\times_{U}U_{\beta}$ and the unlabeled maps are
  the structural isomorphisms associated with pseudofunctoriality of $E$.
\end{itemize}
and which acts on arrows by sending $(g,(g_{\alpha}))$ to $
(g|_{V},(g_{\alpha}))$.

\subsection{Local fibrations}

We are now in a position to describe the maps which will be the
fibrations in our fibration structure.

\begin{definition}\label{def:local_fib}
  A map $p\colon E\to B$ is a \myemph{local fibration} if and only if, for
  every $U$ and cover $\SS$ of $U$, the map
  \begin{align*}
    \Phi_{\SS}\colon E(U) \to\descent{p}{\SS}
  \end{align*}
  described in Section \ref{section:descent_data} above is a weak
  equivalence.
\end{definition}
\begin{example}\label{example:stacks}
  When $p$ is the canonical map $F\to 1$, $\descent{p}{(U_{\alpha})}$
  is the pseudolimit from (\ref{eq:is_stack}) and this map is a local
  fibration if and only if $F$ is a stack.
\end{example}
\begin{example}\label{example:prestack_fibration}
  Let $\mathbf{2}$ be the category with two objects, one
  non-identity arrow and one connected component.  Then, for
  $A\colon \CC^{\myop}\to\cat$ a pseudofunctor, $[\mathbf{2},A]$ denotes
  the cotensor with $\mathbf{2}$.  I.e.,
  $[\mathbf{2},A](U)=A(U)^{\mathbf{2}}$.  $A$ is a prestack if and
  only if the induced map $\langle\partial_{0},\partial_{1}\rangle\colon [\mathbf{2},A]\to A\times
  A$ is a local fibration.
\end{example}
Notice that the map $\Phi_{\SS}$ is always faithful and that we have
the following characterization of local fibrations between prestacks:
\begin{lemma}
  If $E$ and $B$ are prestacks, then $p\colon E\to B$ is a local fibration
  if and only if, for each $U$ and cover $\SS$, $\Phi_{\SS}$ is
  essentially surjective on objects.
  \begin{proof}
    Given a map $(f,f_{\alpha})\colon \Phi_{\SS}(e)\to\Phi_{\SS}(e')$ in
    $\descent{p}{\SS}$ it follows from the fact that $E$ is a prestack
    that the $f_{\alpha}$ possess a unique amalgamation
    $g\colon e\to e'$.  Since $B$ is a prestack we may test locally
    to see that $p(g)=f$.
  \end{proof}
\end{lemma}

\section{The fibration structure on prestacks}\label{section:model_structure}

We will now describe the fibration structure on $\prestacks(\CC)$ for
a site $(\CC,J)$ such that the topology $J$ is precanonical.  We begin by defining what we will call \emph{local weak equivalences}
(this definition can be found in \cite{Porter:M} and similar
definitions appear throughout the literature on stacks and homotopy theory).
\begin{definition}
  A map $h\colon A\to B$ in $\prestacks(\CC)$ is said to be
  \myemph{locally essentially surjective on objects} if and only if
  for any $U$ and $b\in B(U)$ there exists a cover
  \mbox{$\SS=(f_{\alpha}\colon U_{\alpha}\to U)$} of $U$ together with, for each $\alpha$, an
  element $\tilde{b}_{\alpha}\in A(U_{\alpha})$ and an isomorphism
  $\psi_{\alpha}\colon  h(\tilde{b}_{\alpha})\to b|_{\alpha}$.
\end{definition}
\begin{definition}
  A map $h\colon A\to B$ in $\prestacks(\CC)$ is a
  \myemph{local weak equivalence} if it is full, faithful and locally
  essentially surjective on objects.
\end{definition}
Here
being full and faithful means being \emph{pointwise} full and
faithful.

The remainder of this section is devoted to giving a proof of the
following result:
\begin{theorem}\label{theorem:prestacks_fs}
  There is a fibration structure on $\prestacks(\CC)$ with fibrations
  the local fibrations and weak equivalences the local weak equivalences.
\end{theorem}
Consequently, the fibrant objects in this case are precisely the
stacks.
\begin{corollary}\label{theorem:stacks}
  There is an equivalence of 2-categories $\stacks(\CC)\myequiv\hocat\bigl(\prestacks(\CC)\bigr)$.
\end{corollary}
Throughout the remainder of this section we denote by $\frak{W}$ the
class of local weak equivalences and by $\frak{F}$ the class of maps $p$ such that
$\frak{W}\pitchfork p$.

\subsection{Three-for-two}

We will now show that the local weak equivalences satisfy the
three-for-two condition:
\begin{proposition}
  \label{prop:three_for_two}
  Given a diagram
  \begin{align*}
    \xy
    {\ar^{f}(0,0)*+{A};(24,0)*+{C}};
    {\ar_{g}(0,0)*+{A};(12,-14)*+{B}};
    {\ar_{h}(12,-14)*+{B};(24,0)*+{C}};
    {\ar@{=>}^{\gamma}(12,-3)*+{};(12,-8)*+{}};
    \endxy
  \end{align*}
  with $\gamma$ an isomorphism, if any two of $f,g$ and $h$ are local weak
  equivalences, then so is the third.
  \begin{proof}
    If $h$ and $g$ are local weak equivalences, then it is trivial to
    verify that $f$ is also a local weak equivalence.

    When $f$ and $g$ are local weak equivalences it is easily seen
    that $h$ is locally essentially surjective since $f$ is.  To see
    that $h$ is full, suppose given a map $j\colon h(x)\to h(y)$ in $C(U)$.
    Because $g$ is locally essentially surjective on objects we can find
    a cover $\SS$ of $U$ and isomorphisms
    $\varphi_{\alpha}\colon g(a_{\alpha})\to x|_{\alpha}$ and
    $\psi_{\alpha}\colon g(b_{\alpha})\to y|_{\alpha}$ in $B(U_{\alpha})$
    for each $U_{\alpha}\to U$ in the cover.  We can then construct
    composites
    \begin{align*}
      \xy
      {\ar_{\gamma_{a_{\alpha}}}(0,0)*+{f(a_{\alpha})};(0,-15)*+{h(g(a_{\alpha}))}};
      {\ar^{h(\varphi_{\alpha})}(0,-15)*+{h(g(a_{\alpha}))};(20,-15)*+{h(x|_{\alpha})}};
      {\ar(20,-15)*+{h(x|_{\alpha})};(40,-15)*+{h(x)|_{\alpha}}};
      {\ar^{j|_{\alpha}}(40,-15)*+{h(x)|_{\alpha}};(60,-15)*+{h(y)|_{\alpha}}};
      {\ar(60,-15)*+{h(y)|_{\alpha}};(80,-15)*+{h(y|_{\alpha})}};
      {\ar^{h(\psi^{-1}_{\alpha})}(80,-15)*+{h(y|_{\alpha})};(100,-15)*+{h(g(b_{\alpha}))}};
      {\ar_{\gamma^{-1}_{b_{\alpha}}}(100,-15)*+{h(g(b_{\alpha}))};(100,0)*+{f(b_{\alpha})}};
      \endxy
    \end{align*}
    where the unlabelled arrows are the coherence isomorphisms
    associated to $h$.  Since $f$ is full and faithful there exists
    a canonical lift $u_{\alpha}\colon a_{\alpha}\to b_{\alpha}$ in
    $A(U_{\alpha})$ for each $U_{\alpha}\to U$ in the cover $\SS$.
    Using these lifts we similarly obtain maps
    $v_{\alpha}\colon x|_{\alpha}\to y|_{\alpha}$ defined as
    $\psi_{\alpha}\circ g(u_{\alpha})\circ \varphi_{\alpha}^{-1}$.
    These constitute a matching family for the presheaf $B(x,y)$.  To
    see this it suffices to show that, for each $U_{\alpha}\to U$ and
    $U_{\beta}\to U$ in $\SS$, the diagram
    \begin{align*}
      \xy
      {\ar(0,0)*+{g(a_{\alpha}|_{\alpha\beta})};(0,15)*+{g(a_{\alpha})|_{\alpha\beta}}};
      {\ar^{\varphi_{\alpha}|_{\alpha\beta}}(0,15)*+{g(a_{\alpha})|_{\alpha\beta}};(20,15)*+{x|_{\alpha}|_{\alpha\beta}}};
      {\ar(20,15)*+{x|_{\alpha}|_{\alpha\beta}};(35,15)*+{x|_{\alpha\beta}}};
      {\ar(35,15)*+{x|_{\alpha\beta}};(50,15)*+{x|_{\beta}|_{\alpha\beta}}};
      {\ar^{\varphi^{-1}_{\beta}|_{\alpha\beta}}(50,15)*+{x|_{\beta}|_{\alpha\beta}};(75,15)*+{g(a_{\beta})|_{\alpha\beta}}};
      {\ar^{g(u_{\beta})|_{\alpha\beta}}(75,15)*+{g(a_{\beta})|_{\alpha\beta}};(100,15)*+{g(b_{\beta})|_{\alpha\beta}}};
      {\ar(100,15)*+{g(b_{\beta})|_{\alpha\beta}};(100,0)*+{g(b|_{\beta}|_{\alpha\beta})}};
      {\ar(0,0)*+{g(a_{\alpha}|_{\alpha\beta})};(0,-15)*+{g(a_{\alpha})|_{\alpha\beta}}};
      {\ar_{g(u_{\alpha})|_{\alpha\beta}}(0,-15)*+{g(a_{\alpha})|_{\alpha\beta}};(25,-15)*+{g(b_{\alpha})|_{\alpha\beta}}};
      {\ar_{\psi_{\alpha}}(25,-15)*+{g(b_{\alpha})|_{\alpha\beta}};(50,-15)*+{y|_{\alpha}|_{\alpha\beta}}};
      {\ar(50,-15)*+{y|_{\alpha}|_{\alpha\beta}};(65,-15)*+{y|_{\alpha\beta}}};
      {\ar(65,-15)*+{y|_{\alpha\beta}};(80,-15)*+{y|_{\beta}|_{\alpha\beta}}};
      {\ar_{\psi^{-1}_{\beta}}(80,-15)*+{y|_{\beta}|_{\alpha\beta}};(100,-15)*+{g(b|_{\beta})|_{\alpha\beta}}};
      {\ar(100,-15)*+{g(b|_{\beta})|_{\alpha\beta}};(100,0)*+{g(b|_{\beta}|_{\alpha\beta})}};
      \endxy
    \end{align*}
    commutes, where the unnamed arrows are the evident coherence
    isomorphisms.  Since $g$ is full and faithful both ways
    around this diagram induces canonical lifts
    \mbox{$\xi,\zeta\colon a_{\alpha}|_{\alpha\beta}\to
    b_{\beta}|_{\alpha\beta}$}.  It suffices by faithfulness of $f$ to
    show that $f(\xi)=f(\zeta)$, which holds by a straightforward
    diagram chase.  Since the $v_{\alpha}$ are a matching family it
    follows from the fact that $B$ is a prestack that there exists a
    canonical amalgamation $v\colon x\to y$ in $B(U)$.  This map clearly has
    the property that $h(v)=j$, as required.

    To see that $h$ is faithful one uses roughly the same kind of
    approach.  Given $j,k\colon x\to y$ in $B(U)$ with $h(j)=k(j)$ we use
    local essential surjectivity of $g$ to obtain a cover $\SS$ and
    isomorphisms $g(a_{\alpha})\iso x|_{\alpha}$ and
    $g(b_{\alpha})\iso y|_{\alpha}$.  Conjugation of $j|_{\alpha}$ and
    $k|_{\alpha}$ by these isomorphisms gives two families of maps
    $g(a_{\alpha})\to g(b_{\alpha})$ and since $g$ is full and
    faithful these induce canonical lifts
    $u_{\alpha},v_{\alpha}\colon a_{\alpha}\to b_{\alpha}$ in
    $A(U_{\alpha})$.  Using the fact that $h(j)=k(j)$ we can then show
    that $f(u_{\alpha})=f(v_{\alpha})$ so that
    $u_{\alpha}=v_{\alpha}$.  It then follows by the fact that $B$ is
    a prestack that $j=k$.

    The proof that $g$ is a local weak equivalence when $f$ and $h$
    are is similar and is left to the reader.
  \end{proof}
\end{proposition}

\subsection{Characterization of the fibrations}

We now turn to providing a characterization of the fibrations
$\frak{F}$.  This result is analogous to an earlier result of Joyal
and Tierney \cite{Joyal:SSCS} in which they characterize stacks as
weakly fibrant objects.  The differences between our result and theirs
are as follows.  First, they consider a Grothendieck topos $\EE$ with
the canonical topology and they characterize those groupoids
$G$ in $\EE$ such that the externalization $\EE(-,G)$ is a stack.  
In our case, the site is an arbitrary precanonical site and our
prestacks are fibered in categories rather than groupoids.  In the
setting of \emph{ibid} it is not necessary to consider prestacks and
it is not necessary to make use of the axiom of choice.  Because we
work in a more general setting we must restrict first to prestacks and
we also appeal to the axiom of choice. Finally, the characterization
we give is of local fibrations in general and not just stacks.
\begin{lemma}
  \label{lemma:w_e_descent}
  For $i\colon A\to C$ in $\frak{W}$ and $U$ in $\CC$, every
  object $c$ of $C(U)$ determines a cover $\SS$ and an object of
  $\descent{i}{\SS}$ which projects via $\pi\colon \descent{i}{\SS}\to C(U)$
  onto $c$.
  \begin{proof}
    Let an object $c$ of $C(U)$ be given.  Because $i$ is locally
    essentially surjective on objects there exists a family of
    isomorphisms $\psi_{\alpha}\colon i(\tilde{c}_{\alpha})\to c|_{\alpha}$.
    We may form the composites
    \begin{align*}
      \xy
      {\ar(0,0)*+{i(\tilde{c}_{\beta}|_{\alpha\beta})};(22,0)*+{i(\tilde{c}_{\beta})|_{\alpha\beta}}};
      {\ar^{\psi_{\beta}|_{\alpha\beta}}(22,0)*+{i(\tilde{c}_{\beta})|_{\alpha\beta}};(42,0)*+{c|_{\beta}|_{\alpha\beta}}};
      {\ar^{\sigma_{\alpha\beta}(c)}(42,0)*+{c|_{\beta}|_{\alpha\beta}};(62,0)*+{c|_{\alpha}|_{\alpha\beta}}};
      {\ar^{\psi^{-1}_{\alpha}|_{\alpha\beta}}(62,0)*+{c|_{\alpha}|_{\alpha\beta}};(82,0)*+{i(\tilde{c}_{\alpha})|_{\alpha\beta}}};
      {\ar(82,0)*+{i(\tilde{c}_{\alpha})|_{\alpha\beta}};(105,0)*+{i(\tilde{c}_{\alpha}|_{\alpha\beta})}};
      \endxy
    \end{align*}
    where the unlabelled arrows are induced by the coherence 2-cell associated to the
    pseudonatural transformation $i$.
    Since $i$ is full and faithful these possess canonical
    invertible lifts
    \mbox{$\vartheta_{\alpha\beta}\colon \tilde{c}_{\beta}|_{\alpha\beta}\to\tilde{c}_{\alpha}|_{\alpha\beta}$}
    in $A(U_{\alpha\beta})$.  It is routine to verify that
    $(c,(\tilde{c}_{\alpha}),(\psi_{\alpha}),(\vartheta_{\alpha\beta}))$
    is an object of $\descent{i}{\SS}$.
  \end{proof}
\end{lemma}

\begin{lemma}
  \label{lemma:local_fib_lifting}
  If $p\colon E\to B$ is a local fibration, then
  $\frak{W}\pitchfork p$.
  \begin{proof}
    Suppose $p\colon E\to B$ is a local fibration and let a diagram of the
    form
    \begin{align*}
      \xy
      {\ar^{h}(0,20)*+{A};(20,20)*+{E}};
      {\ar_{i}(0,20)*+{A};(0,0)*+{C}};
      {\ar^{p}(20,20)*+{E};(20,0)*+{B}};
      {\ar_{k}(0,0)*+{C};(20,0)*+{B}};
      {\ar@{=>}^{\gamma}(10,13)*+{};(10,7)*+{}};
      \endxy
    \end{align*}
    be given with $i\colon A\to C$ in $\frak{W}$.  Given an
    object $c$ of $C(U)$ we may choose, by Lemma \ref{lemma:w_e_descent},
    a cover $\SS$ together with descent data
    $(c,(a_{\alpha}),(\psi_{\alpha}),(\vartheta_{\alpha\beta}))$ in
    $\descent{i}{\SS}$.  By Lemma \ref{lemma:functorial_squares} this
    gives descent data
    $\hat{c}:=(k(c),(h(a_{\alpha})),(\hat{\psi}_{\alpha}),(\hat{\vartheta}_{\alpha\beta}))$
    in $\descent{p}{\SS}$.  Thus, we choose $l(c)$ to be an
    amalgamation of this descent data.

    Given $f\colon c\to d$ in $C(U)$ assume that $\SS$ and $\RR$
    are the covers chosen in the definition of $l(c)$ and $l(d)$
    and assume that
    $(c,(a_{\alpha}),(\psi_{\alpha}),(\vartheta_{\alpha\beta}))$ and
    $(d,(b_{\alpha}),(\varphi_{\alpha}),(\omega_{\alpha\beta}))$ are
    the descent data chosen in the definition of $l(c)$ and $l(d)$, respectively.
    Let $\WW$ be the common refinement $\SS\cap\RR$ of $\SS$ and $\RR$ and observe
    that, for $U_{\alpha}\to U$ in $\WW$, we have isomorphisms
    $\chi\colon \Phi_{\SS}(l(c))\iso\hat{c}$ and
    $\mu\colon \Phi_{\RR}(l(d))\iso\hat{d}$.  We also have
    \begin{align*}
      \xy
      {\ar^{\psi_{\alpha}}(0,0)*+{i(a_{\alpha})};(20,0)*+{c|_{\alpha}}};
      {\ar^{f|_{\alpha}}(20,0)*+{c|_{\alpha}};(40,0)*+{d|_{\alpha}}};
      {\ar^{\varphi^{-1}_{\alpha}}(40,0)*+{d|_{\alpha}};(60,0)*+{i(b_{\alpha}).}};
      \endxy
    \end{align*}
    As such, since $i\colon A\to C$ is full and
    faithful, there exists a canonical map $f_{\alpha}\colon a_{\alpha}\to
    b_{\alpha}$ which is mapped by $i$ onto this composite.  This
    gives a map of descent data
    \begin{align*}
      (f,f_{\alpha})\colon (c,(a_{\alpha}),(\psi_{\alpha}),(\vartheta_{\alpha\beta}))|_{\WW}\to(d,(b_{\alpha}),(\varphi_{\alpha}),(\omega_{\alpha\beta}))|_{\WW}
    \end{align*}
    and by Lemma \ref{lemma:functorial_squares}, we have that
    $(k(f),(h(f_{\alpha})))\colon \hat{c}|_{\WW}\to\hat{d}|_{\WW}$ in
    $\descent{p}{\WW}$.  Therefore we may form the composite
    \begin{align*}
      \xy
      {\ar^(.6){\chi|_{\WW}}(0,0)*+{\Phi_{\SS}(l(c))|_{\WW}};(25,0)*+{\hat{c}|_{\WW}}};
      {\ar^{(k(f),(h(f_{\alpha})))}(25,0)*+{\hat{c}|_{\WW}};(50,0)*+{\hat{d}|_{\WW}}};
      {\ar^(.4){\mu^{-1}|_{\WW}}(50,0)*+{\hat{d}|_{\WW}};(75,0)*+{\Phi_{\VV}(l(d))|_{\WW}}};
      \endxy
    \end{align*}
    which gives us a family of maps $l(c)|_{\alpha}\to
    l(d)|_{\alpha}$ for $U_{\alpha}\to U$ in $\WW$.  This family
    constitutes a matching family for $E(l(c),l(d))$ and since $E$ is
    a prestack there exists a canonical amalgamation $l(f)\colon l(c)\to
    l(d)$.  Functoriality follows from the uniqueness of amalgamations.

    We now construct the natural isomorphisms $\lambda\colon h\iso l\circ i$
    and $\rho\colon p\circ l\iso k$.  First, for $\lambda$, assume given an
    object $u$ of $AU$.  Assume that $\SS$ is the cover of $U$
    and $(iu,(a_{\alpha}),(\psi_{\alpha}),(\vartheta_{\alpha\beta}))$
    is the descent data chosen in the construction of $l(iu)$.  Then
    $(1_{iu},(\psi^{-1}_{\alpha}))$ is an isomorphism in $\descent{i}{\SS}$ from
    $\Phi_{\SS}(u)$ to
    $(iu,(a_{\alpha}),(\psi_{\alpha}),(\vartheta_{\alpha\beta}))$.  As
    such, we may form the following composite
    \begin{align*}
      \xy
      {\ar^(.4){\hat{\gamma}(a)}(0,0)*+{\Phi_{\SS}(ha)};(30,0)*+{(h,k,\gamma)_{*}(\Phi_{\SS}(a))}};
      {\ar^(.6){(h,k,\gamma)_{*}(1_{iu},(\psi^{-1}_{\alpha}))}(30,0)*+{(h,k,\gamma)_{*}(\Phi_{\SS}(a))};(70,0)*+{\widehat{ia}}};
      {\ar(70,0)*+{\widehat{ia}};(90,0)*+{\Phi_{\SS}(l(ia))}};
      \endxy
    \end{align*}
    in $\descent{\SS}{p}$, where $\hat{\gamma}$ is as in the
    discussion of $(h,k,\gamma)_{*}$ from Section
    \ref{section:descent_data} and the unnamed map is the isomorphism
    associated to the definition of $l(ia)$.  Because $\Phi_{U}$ is full
    and faithful this gives a canonical isomorphism
    $\lambda(a)\colon h(a)\to l(ia)$ with $\Phi_{\SS}(\lambda(a))$ the
    composite above.  Naturality of $\lambda$ follows from
    faithfulness of the $\Phi_{\SS}$ together with the definition of
    the action of $l$ on arrows.  Next, we define $\rho(c)\colon p(l(c))\to k(c)$ to be the first component of
    the isomorphism \mbox{$\Phi_{\SS}(l(c))\iso \hat{c}$} of descent
    data associated to the definition of $l(c)$.  This is natural by
    definition of $l$.  Finally, it is immediate from the definitions that $\gamma$ can be
    recovered by composing the pasting diagram obtained from $\lambda$
    and $\rho$.
  \end{proof}
\end{lemma}
\begin{theorem}\label{theorem:fib_char}
  For a map $p\colon E\to B$ the following are equivalent:
  \begin{enumerate}
  \item $p$ is a local fibration.
  \item $p$ is in $\frak{F}$.
  \end{enumerate}
  \begin{proof}
    By Lemma \ref{lemma:local_fib_lifting} it suffices to prove that
    if $p\colon E\to B$ is in $\frak{F}$, then it is a local
    fibration.  To this end, let $U$ together with a cover $\SS$
    be given.  Assume given descent data
    $(b,(e_{\alpha}),(\psi_{\alpha}),(\vartheta_{\alpha\beta}))$ in
    $\descent{p}{\SS}$.  Then we have a square
    \begin{align*}
      \xy
      {\ar^{e}(0,20)*+{\hat{\SS}};(20,20)*+{E}};
      {\ar_{i}(0,20)*+{\hat{\SS}};(0,0)*+{yU}};
      {\ar_{b}(0,0)*+{yU};(20,0)*+{B}};
      {\ar^{p}(20,20)*+{E};(20,0)*+{B}};
      {\ar@{=>}^{\psi}(10,13)*+{};(10,7)*+{}};
      \endxy
    \end{align*}
    where $yU$ is the representable functor and $\hat{\SS}$ is the
    subfunctor of $yU$ induced by the cover
    $\SS$ (note that both of these are prestacks).  Also, $e$ is the
    pseudonatural transformation representing
    the family $(e_{\alpha})$ with coherence isomorphisms constructed
    using the $\vartheta_{\alpha\beta}$.  Similarly, $b$ is the
    pseudonatural transformation representing $b$.  Finally, $\psi$ is
    the modification with component at $U_{\alpha}\to
    U$ in $\SS$ given by $\psi_{\alpha}$.
    
    Notice that $i$ is a local weak equivalence so that, since
    $\frak{W}\pitchfork p$, it follows that there exists
    a lift $l\colon yU\to E$ together with isomorphisms
    $\lambda\colon e\iso l\circ i$ and $\rho\colon p\circ l\iso b$ such that the
    square above can be recovered from these.  I.e., we have $l$ an
    object of $EU$ together with $\rho_{V}\colon p(l)|_{V}\iso b|_{V}$ for
    every $V\to U$
    and $\lambda_{\alpha}\colon e_{\alpha}\iso l|_{\alpha}$ for each
    $U_{\alpha}\to U$ in the cover.  It is then routine to verify that $l$
    is an amalgamation of our descent data.
    \end{proof}
\end{theorem}
\begin{corollary}
  For any $F$, the canonical map $F\to 1$ is in $\frak{F}$ if and only if $F$ is a stack.
\end{corollary}
\begin{corollary}
  If $p\colon E\to B$ is in $\frak{F}\cap\frak{W}$, then $p$ is an
  equivalence (i.e., there exists a $p'\colon B\to E$ together with
  invertible $\eta\colon 1_{B}\to p\circ p'$ and $\epsilon\colon p'\circ p\to 1_{E}$).
\end{corollary}
\begin{corollary}
  \label{corollary:wild_fact}
  If $p\colon E\to B$ is in $\frak{F}\cap\frak{W}$ and $i\colon A\to C$ is any
  map, then $i\pitchfork p$.
\end{corollary}
\begin{corollary}
  Theorem \ref{theorem:fib_char} is equivalent to the Axiom of Choice.
  \begin{proof}
    Consider the case where our site consists of the lattice $\OO(\emptyset)$
    of open subsets of the empty set with its
    canonical topology and the notion of covering family is given by
    the usual topological notion of covering family.  In this case we are working directly in
    $\cat$ and we can easily prove that every object is locally fibrant.
    Using this it is possible to construct pseudo-inverses of weak
    categorical equivalences.  Therefore the Axiom of Choice holds.
  \end{proof}
\end{corollary}

\subsection{Factorization and isocomma objects}\label{section:factorizations}

We will now describe the factorizations in $\prestacks(\CC)$.  To a
map $f\colon A\to B$ we associate a prestack
$\mypath(f)$ by letting $\mypath(f)(U)$ be the category with
\begin{description}
\item[Objects] Tuples consisting of a cover $\SS$ and an object
  $(b,(e_{\alpha}),(\psi_{\alpha}),(\vartheta_{\alpha\beta}))$ of $\descent{f}{\SS}$.
\item[Arrows] An arrow
  $(\SS,b,(e_{\alpha}),(\psi_{\alpha}),(\vartheta_{\alpha\beta}))\to(\VV,b',(e'_{\alpha}),(\psi'_{\alpha}),(\vartheta'_{\alpha\beta}))$
  is an equivalence class of data given by a common refinement $\WW$ of $\SS$ and $\VV$ together with a map
  \begin{align*}
    \bigl(b,(e_{\alpha}),(\psi_{\alpha}),(\vartheta_{\alpha\beta})\bigr)|_{\WW}\to\bigl(b',(e'_{\alpha}),(\psi'_{\alpha}),(\vartheta'_{\alpha\beta})\bigr)|_{\WW}
  \end{align*}
  in $\descent{f}{\WW}$.  We identify $(\WW,g,(g_{\alpha}))$ and $(\WW',g',(g_{\alpha}'))$ when
  there exists a common refinement of $\WW$ and $\WW'$ on which the
  maps of descent data agree.
\end{description}
Note that $g=g'$ when $(\WW,g,(g_{\alpha}))$ and $(\WW',g',(g_{\alpha}'))$
are identified in $\mypath{f}(U)$. 

There is, for $g\colon V\to U$, an obvious restriction map
$\mypath(f)(U)\to\mypath(f)(V)$ which acts by pullback on both covers
and descent data.  This makes $\mypath(f)$ into a pseudofunctor
$\CC^{\myop}\to\cat$.  We observe that we have the following lemma,
the proof of which is straightforward:
\begin{lemma}
  \label{lemma:graph_prestack}
  If $A$ and $B$ are prestacks and $f\colon A\to B$, then $\mypath(f)$ is a prestack.
\end{lemma}
There is a projection map $\mypath(f)\to B$ which sends
$(\SS,b,(e_{\alpha}),(\psi_{\alpha}),(\vartheta_{\alpha\beta}))$ to
$b$ and sends an arrow $[\WW,g,(g_{\alpha})]$ to $g$.  We define
$Q\colon\prestacks(\CC)^{\to}\to\prestacks(\CC)^{\to}$
by letting $Q(f)$, for $f\colon A\to B$ in $\prestacks(\CC)$, be the
projection $\mypath(f)\to B$.

For the pseudonatural transformation $\eta\colon
1_{\prestacks(\CC)^{\to}}\to Q$, note that there is a map
$A\to\mypath(f)$, which we denote by $\eta_{f}$, that sends an $a$ in
$A(U)$ to $(M_{U},f(a),(a|_{\alpha}),(\psi_{\alpha}),(\vartheta_{\alpha\beta}))$
where $M_{U}$ denotes the maximal sieve on $U$, the $\psi_{\alpha}$
are the coherence isomorphisms associated to
$f$, and the $\vartheta_{\alpha\beta}$ are the coherence isomorphisms
obtained from the structure of $A$ as a pseudofunctor.  It is
straightforward to verify that $Q(f)\circ\eta_{f}=f$ and this equation
determines the rest of the data of the pseudonatural transformation $\eta$.
\begin{lemma}\label{lemma:factorization_1}
  For $f\colon A\to B$, in the factorization
    \begin{align*}
    \xy
    {\ar^{f}(0,15)*+{A};(30,15)*+{B}};
    {\ar_{\eta_{f}}(0,15)*+{A};(15,0)*+{\mypath(f)}};
    {\ar_{Q(f)}(15,0)*+{\mypath(f)};(30,15)*+{B}};
    \endxy
  \end{align*}
  $Q(f)$ is a local fibration and $\eta_{f}$ is in $\frak{W}$.
  \begin{proof}
      It is trivial that $\eta_{f}$ is in $\frak{W}$.  To see
      that $Q(f)$ is a local fibration let a cover $\WW=(h^{\gamma}\colon U^{\gamma}\to U)_{\gamma}$ of $U$ be
      given together with an object
      \begin{align*}
        b,\bigl(\SS^{\gamma},b^{\gamma},(e^{\gamma}_{\alpha}),(\psi^{\gamma}_{\alpha}),(\vartheta^{\gamma}_{\alpha\beta})\bigr),(\varphi_{\gamma}),(\Theta_{\gamma\delta})
      \end{align*}
      of $\descent{Q(f)}{\WW}$ where
      $\SS^{\gamma}=(h^{\gamma}_{\alpha}\colon U^{\gamma}_{\alpha}\to U^{\gamma})_{\alpha}$.  Here
      $\Theta_{\gamma\delta}=(h_{\gamma\delta},(h^{\alpha}_{\gamma\delta}))$
      is an isomorphism
      \begin{align*}
        \xy
        {\ar^{\iso}(0,0)*+{(b_{\delta},(e^{\delta}_{\alpha}),(\psi^{\delta}_{\alpha}),(\vartheta^{\delta}_{\alpha\beta}))|_{\gamma\delta}};(45,0)*+{(b_{\gamma},(e^{\gamma}_{\alpha}),(\psi^{\gamma}_{\alpha}),(\vartheta^{\gamma}_{\alpha\beta}))|_{\gamma\delta}}};
        \endxy
      \end{align*}
      of descent data in $\descent{f}{\SS^{\gamma\delta}}$.  Define a
      new cover $\bar{\WW}$ of $U$ as the cover consisting of the maps
      of the form $h^{\gamma}\circ h^{\gamma}_{\alpha}\colon U^{\gamma}_{\alpha}\to
      U^{\gamma}\to U$ for $h^{\gamma}$ in $\WW$ and
      $h^{\gamma}_{\alpha}$ in $\UU^{\gamma}$.  We then have an
      object 
      \begin{align*}
        (\bar{\WW},b,(e^{\gamma}_{\alpha}),(\varphi^{\gamma}_{\alpha}),(\chi^{\gamma\delta}_{\alpha\beta}))
      \end{align*}
      of $\descent{f}{\SS}$ where $\varphi^{\gamma}_{\alpha}$ is the
      composite
      \begin{align*}
        \xy
        {\ar^{\psi^{\gamma}_{\alpha}}(0,0)*+{f(e_{\alpha}^{\gamma})};(20,0)*+{b^{\gamma}|_{U^{\gamma}_{\alpha}}}};
        {\ar^{\varphi_{\gamma}|_{U^{\gamma}_{\alpha}}}(20,0)*+{b^{\gamma}|_{U^{\gamma}_{\alpha}}};(40,0)*+{b|_{U^{\gamma}_{\alpha}}}};
        \endxy
      \end{align*}
      and $\chi^{\gamma\delta}_{\alpha\beta}$ is the composite
      \begin{align*}
        \xy
        {\ar^{h^{\beta}_{\gamma\delta}}(0,0)*+{e^{\delta}_{\beta}|_{U^{\gamma}_{\alpha}\cap
              U^{\delta}_{\beta}}};(30,0)*+{e^{\gamma}_{\beta}|_{U^{\gamma}_{\alpha}\cap
              U^{\delta}_{\beta}}}};
        {\ar^{\vartheta^{\gamma}_{\alpha\beta}}(30,0)*+{e^{\gamma}_{\beta}|_{U^{\gamma}_{\alpha}\cap
              U^{\delta}_{\beta}}};(60,0)*+{e^{\gamma}_{\alpha}|_{U^{\gamma}_{\alpha}\cap U^{\delta}_{\beta}}.}};
        \endxy
      \end{align*}
      With these definitions, it is a (quite) lengthy but straightforward
      verification that we have described the amalgamation of the
      descent data.
 \end{proof}
\end{lemma}
\begin{example}
  When $f$ is the canonical map $A\to 1$ we see that $\mypath(f)$ is the
  associated stack $\astack(A)$ of $A$
  (cf.~\cite{Moerdijk:ILSG,Porter:M} for more on the associated stack).
\end{example}
We note that when $A$ is a stack it is possible to factor $f$ in a
more straightforward way using isocomma objects.  
\begin{definition}
  Given maps $f\colon A\to B$ and $g\colon C\to B$ in
  $\prestacks(\CC)$, the \myemph{isocomma object $(f,g)$} is the
  pseudofunctor given at $U$ by the category $(f,g)(U)$ with 
  \begin{description}
  \item[Objects] Tuples consisting of objects $a$ and $c$ of $A(U)$
    and $C(U)$, respectively, and an isomorphism $\xi\colon f(a)\iso g(c)$.
  \item[Arrows] An arrow $(a,c,\xi)\to(b,d,\zeta)$ is given by maps
    $i\colon a\to b$ and $j\colon c\to d$ such that the diagram
    \begin{align*}
      \xy
      {\ar^{f(i)}(0,15)*+{f(a)};(20,15)*+{f(b)}};
      {\ar_{\xi}(0,15)*+{f(a)};(0,0)*+{g(c)}};
      {\ar_{g(j)}(0,0)*+{g(c)};(20,0)*+{g(d)}};
      {\ar^{\zeta}(20,15)*+{f(b)};(20,0)*+{g(d)}};
      \endxy
    \end{align*}
    commutes.
  \end{description}
  The action of $(f,g)$ on arrows is simply by restriction of all of
  the aforementioned data.
\end{definition}
There is an invertible 2-cell $\chi$ as indicated in the following diagram:
\begin{align*}
  \xy
  {\ar(0,15)*+{(f,g)};(20,15)*+{C}};
  {\ar(0,15)*+{(f,g)};(0,0)*+{A}};
  {\ar^{g}(20,15)*+{C};(20,0)*+{B}};
  {\ar_{f}(0,0)*+{A};(20,0)*+{B}};
  {\ar@{=>}^{\chi}(10,10)*+{};(10,4)*+{}};
  \endxy
\end{align*}
where the unnamed arrows are the obvious projections.  Here $\chi$
projects $(a,c,\xi)\mapsto\xi$.  The universal property of $(f,g)$ is
that for any other diagram
\begin{align*}
  \xy
  {\ar(0,15)*+{Z};(20,15)*+{C}};
  {\ar(0,15)*+{Z};(0,0)*+{A}};
  {\ar^{g}(20,15)*+{C};(20,0)*+{B}};
  {\ar_{f}(0,0)*+{A};(20,0)*+{B}};
  {\ar@{=>}^{\chi'}(10,10)*+{};(10,4)*+{}};
  \endxy
\end{align*}
with $\chi'$ invertible, there exists a canonical map $z\colon Z\to
(f,g)$ such that the diagram
\begin{align*}
  \xy
  {\ar^{z}(15,10)*+{Z};(15,0)*+{(f,g)}};
  {\ar@/_1em/(15,10)*+{Z};(0,0)*+{A}};
  {\ar(15,0)*+{(f,g)};(0,0)*+{A}};
  {\ar@/^1em/(15,10)*+{Z};(30,0)*+{C}};
  {\ar(15,0)*+{(f,g)};(30,0)*+{C}};
  \endxy
\end{align*}
commutes and such that $\chi'=\chi\circ z$.  It is straightforward to
show that $(f,g)$ is a prestack when $A$ and $C$ are.

Now, the universal property gives us a map $i\colon A\to (f,1_{B})$
such that 
\begin{align*}
  \xy
  {\ar^{f}(0,10)*+{A};(20,10)*+{B}};
  {\ar_{i}(0,10)*+{A};(10,0)*+{(f,1_{B})}};
  {\ar_{p}(10,0)*+{(f,1_{B})};(20,10)*+{B}};
  \endxy
\end{align*}
commutes, where $p$ is the projection.  Here it is clear that this
gives a factorization $f=p\circ i$.  In particular, $i(a)$ is
$(a,f(a),1_{f(a)})$ and it is straightforward to verify that $i$ is in $\frak{W}$.
\begin{lemma}
  When $A$ is a stack, $p\colon (f,1_{B})\to B$ is a local fibration.
  \begin{proof}
    Given descent data $(b,(e_{\alpha}\colon f(a^{\alpha})\iso
    b^{\alpha}),(\psi_{\alpha}),(\vartheta_{\alpha\beta}))$ in
    $\descent{\SS}{p}$, we have that $\vartheta_{\alpha\beta}$ is a
    commutative square
    \begin{align*}
      \xy
      {\ar^{f(\chi_{\alpha\beta})}(0,15)*+{f(a^{\beta})|_{\alpha\beta}};(25,15)*+{f(a^{\alpha})|_{\alpha\beta}}};
      {\ar_{e_{\beta}}(0,15)*+{f(a^{\beta})|_{\alpha\beta}};(0,0)*+{b^{\beta}|_{\alpha\beta}}};
      {\ar^{e_{\alpha}}(25,15)*+{f(a^{\alpha})|_{\alpha\beta}};(25,0)*+{b^{\alpha}|_{\alpha\beta}}};
      {\ar_{\omega_{\alpha\beta}}(0,0)*+{b^{\beta}|_{\alpha\beta}};(25,0)*+{b^{\alpha}|_{\alpha\beta}}};
      \endxy
    \end{align*}
    of isomorphisms.  This gives us descent data
    $((a^{\alpha}),(\chi_{\alpha\beta}))$ for $A$ and $\SS$ and
    since $A$ is a stack there is an amalgamating object $a$ of
    $A(U)$.  For each $\alpha$, we have the isomorphism
    \begin{align*}
      \xy
      {\ar(0,0)*+{f(a)|_{\alpha}};(20,0)*+{f(a|_{\alpha})}};
      {\ar(20,0)*+{f(a|_{\alpha})};(40,0)*+{f(a^{\alpha})}};
      {\ar^{e_{\alpha}}(40,0)*+{f(a^{\alpha})};(60,0)*+{b^{\alpha}}};
      {\ar^{\psi_{\alpha}}(60,0)*+{b^{\alpha}};(80,0)*+{b|_{\alpha}}};
      \endxy
    \end{align*}
    and these are easily seen to constitute a matching family for
    $B(f(a),b)$.  Therefore, since $B$ is a prestack there is a
    canonical amalgam $e\colon f(a)\iso b$.  We define this isomorphism to
    be the object of $(f,B)(U)$ corresponding to our descent data.  It
    is routine to verify that this constitutes a pseudo-inverse to the
    map $(f,B)(U)\to\descent{\SS}{p}$ satisfying the coherence
    conditions from the definition of local fibrations.
  \end{proof}
\end{lemma}
This completes the proof of Theorem \ref{theorem:prestacks_fs}.

\section{Topological, differentiable and algebraic stacks}\label{section:examples}

We will now show that the results of Section
\ref{section:model_structure} can be used to give analogous
characterizations of the 2-categories of topological, differentiable
and algebraic stacks.  These three cases are formal analogues.  The
categories of topological spaces, differentiable manifolds and schemes
all have in common that quotients in them are not well-behaved.  This
gives rise to the situation, familiar from the theory of \'{e}tendues
from \cite{SGA4}, in which one would
like to form a ``generalized quotient'' of a space, manifold or
scheme.  (Indeed, there is an important connection with the theory of
\'{e}tendues as described in \cite{Pronk:ESBF}, but we do
not describe it here.)  Topological, differentiable and algebraic stacks are the
appropriate ``generalized quotients'' of suitable equivalence
relations in each of these situations.  These three cases are formally
analogous in the sense that topological, differentiable and algebraic stacks are by
definition stacks $X$ which appear in a suitable sense as
``quotients'' of topological spaces, differentiable manifolds, or
schemes, respectively.  This formal analogy permits us to give a
single argument (here described in detail for topological stacks)
which will show that each of these 2-categories can be described as
the homotopy 2-category of the corresponding 2-category of prestacks.

\subsection{Topological stacks}

We will briefly recall the definition of topological stacks, which are
essentially the topological version of the \emph{algebraic stacks} of
Deligne and Mumford \cite{Deligne:ISCGG}.  Throughout this section we will be
working with the \myemph{topological site} which consists of a small
category $\topcat$ of sober topological spaces $U,V,\ldots$ equipped
with the \'{e}tale Grothendieck topology.  The \'{e}tale topology is generated by
families $(f_{i}\colon U_{i}\to U)_{i}$ which are said to cover when the
map $\sum_{i}U_{i}\to U$ is an \'{e}tale surjection.
\begin{definition}
  A map $f\colon A\to B$ of prestacks is \myemph{representable} if, for any space $U$
  in $\topcat$ and map $g\colon  yU\to B$, the isocomma object $(f,g)$ is
  representable.
\end{definition}
We now consider pseudofunctors $\pseudo{\topcat}{\groupoids}$ valued
in groupoids.  Throughout this section ``prestack'' means prestack valued
in groupoids and similarly for ``stack''.  Roughly, topological
prestacks are those prestacks which arise as quotients of spaces.
\begin{definition}\label{def:top_prestack}
  A \myemph{topological prestack} is a prestack $A$ such that the
  following conditions are satisfied:
  \begin{enumerate}
  \item The diagonal $\Delta\colon  A\to A\times A$ is representable.
  \item There exists a space $U$ in $\topcat$ and a map $q\colon yU\to A$
    such that, for all spaces $V$ in $\topcat$ and maps $f\colon yV\to A$, the
    map $(f,q)\to V$ is an \'{e}tale surjection.
  \end{enumerate}
\end{definition}
Notice that it makes sense in condition (2) to say that
$(f,q)\to V$ is an \'{e}tale surjection since the domain of
this map is, by condition (1), representable.  We will often refer to
the map $q\colon yU\to A$ from condition (2) as a \myemph{chart for
  $A$}.  Observe that representables are trivially topological
prestacks.  We will henceforth omit explicit mention of the Yoneda
embedding $y$ when no confusion will result.

We denote by $\toppre$ the 2-category of topological prestacks
and we observe that it is an immediate consequence of Lemma
\ref{lemma:local_fib_lifting} that if $p\colon E\to B$ is a local fibration
between topological prestacks, then $\frak{W}\pitchfork p$
where $\frak{W}$ denotes the class of local weak
equivalences in $\toppre$.  We will now consider to what extent
the additional structure of $\prestacks(\topcat)$ restricts to
$\toppre$.
\begin{lemma}
  \label{lemma:top_stack_rep_diag}
  If $f\colon A\to B$ is an equivalence between prestacks and
  $B$ has a representable diagonal, then so does $A$.
  \begin{proof}
    Let maps $v\colon V\to A$ and $w\colon W\to A$ be given.
    Because $B$ has a representable diagonal the isocomma object
    $(f\circ v,f\circ w)$ is a representable $U$.  This gives us the
    following diagram of invertible 2-cells:
    \begin{align*}
      \xy
      {\ar(0,20)*+{U};(20,20)*+{W}};
      {\ar(0,20)*+{U};(0,0)*+{V}};
      {\ar_{v}(0,0)*+{V};(10,0)*+{A}};
      {\ar^{w}(20,20)*+{W};(20,10)*+{A}};
      {\ar^{f}(20,10)*+{A};(20,0)*+{B}};
      {\ar^{f}(10,0)*+{A};(20,0)*+{B}};
      {\ar_{f'}(20,0)*+{B};(30,-10)*+{A}};
      {\ar@/_1.5em/_{1_{A}}(10,0)*+{A};(30,-10)*+{A}};
      {\ar@/^1.5em/^{1_{A}}(20,10)*+{A};(30,-10)*+{A}};
      {(18,-5)*+{\iso}};
      {(25,2)*+{\iso}};
      {(10,10)*+{\iso}};
      \endxy
    \end{align*}
    where $f'$ is a pseudoinverse of $f$.  This is easily seen to exhibit $U$ as $(v,w)$.
  \end{proof}
\end{lemma}

\begin{lemma}
  \label{lemma:top_stack_equiv}
  If $f\colon A\to B$ is an equivalence between prestacks and
  $B$ is a topological prestack, then $A$ is also a topological
  prestack.
  \begin{proof}
    By Lemma \ref{lemma:top_stack_rep_diag} it suffices to prove that
    there exists a space $U$ and an \'{e}tale map $U\to A$.  Because $B$ is a topological
    prestack there exists an \'{e}tale map $e\colon U\to B$.  We
    claim that the map $f'\circ e\colon U\to A$ is \'{e}tale, where
    $f'$ is a pseudoinverse of $f$.  Let another map $v\colon V\to A$
    be given.  Then the isocomma object $(f\circ v,e)$ is a
    representable $W$.  We then obtain the diagram
    \begin{align*}
      \xy
      {\ar(0,15)*+{W};(20,15)*+{U}};
      {\ar(0,15)*+{W};(0,0)*+{V}};
      {\ar_{v}(0,0)*+{V};(10,0)*+{A}};
      {\ar^{f}(10,0)*+{A};(20,0)*+{B}};
      {\ar^{e}(20,15)*+{U};(20,0)*+{B}};
      {\ar_{1_{A}}@/_1em/(10,0)*+{A};(20,-15)*+{A}};
      {\ar^{f'}(20,0)*+{B};(20,-15)*+{A}};
      {(10,7.5)*+{\iso}};
      {(15,-7.5)*+{\iso}};
      \endxy
    \end{align*}
    where the vertical map $W\to V$ is an \'{e}tale surjection.  It is
    straightforward to show that the diagram above exhibits $U$ as
    the isocomma object $(v,f'\circ e)$ so that $f'\circ e$ is
    \'{e}tale.
  \end{proof}
\end{lemma}
Modifying a construction of \cite{Pronk:ESBF}, we associate to each
topological prestack $A$ and chart $e\colon U\to A$ the \'{e}tale
groupoid $G^{e}$ with space of objects $U$
and space of arrows the space representing the isocomma object
$(e,e)$.  In \emph{ibid} it is assumed that $A$ is a
topological stack, but it is in fact sufficient for $A$ to be a
topological prestack.  Also in \emph{ibid} it is shown how to
associate to any \'{e}tale groupoid $G$ a topological stack
$R(G)$. Combining these two procedures, we obtain, for each topological
prestack $A$ and chart $e\colon U\to A$, a topological stack $Q(A,e)$
given by $R(G^{e})$.  In elementary terms, we have
\begin{align*}
  Q(A,e)_{V} & := \textnormal{\textbf{GeomMorph}}\bigl(\sheaves(V),\sheaves(G^{e})\bigr)
\end{align*}
for $V$ a space.  Here the objects are geometric morphisms, arrows are
invertible natural transformations, $\sheaves(V)$ is the ordinary
category of sheaves on the space $V$ and $\sheaves(G^{e})$ is the
category of equivariant sheaves on the groupoid $G^{e}$.  Note that it
is shown in \emph{ibid} that there is a map $i\colon A\to Q(A,e)$ which is a
weak equivalence.

\begin{lemma}\label{lemma:top_fact}
  The associated stack $\astack(A)$ of a topological prestack is a
  topological stack.
  \begin{proof}
    It suffices by Lemma \ref{lemma:top_stack_equiv}, and the fact
    that both $\astack(A)$ and $Q(A,e)$ are both stacks, to construct a
    local weak equivalence $\astack(A)\to Q(A,e)$.  Because the map
    $\eta\colon A\to\astack(A)$ is a local weak equivalence and
    $Q(A,e)$ is a stack there exists a map $\astack(A)\to Q(A,e)$ and
    an invertible 2-cell as indicated in the diagram:
    \begin{align*}
      \xy
      {\ar^(.4){i}(0,15)*+{A};(20,15)*+{Q(A,e).}};
      {\ar@/_1em/_{\eta}(0,15)*+{A};(10,0)*+{\astack(A)}};
      {\ar@/_1em/(10,0)*+{\astack(A)};(20,15)*+{Q(A,e).}};
      {(10,7.5)*+{\iso}};
      \endxy
    \end{align*}
    By the three-for-two property for local weak equivalences it then
    follows that $\astack(A)\to Q(A,e)$ is also a local weak
    equivalence.
  \end{proof}
\end{lemma}
Putting these lemmas together with Theorem \ref{theorem:prestacks_fs}
we have proved the following:
\begin{theorem}\label{theorem:toppre}
  There is a system of fibrant objects on $\toppre$ given by taking
  the local weak equivalences and with fibrant replacement given by
  the associated stack.
\end{theorem}
\begin{corollary}
  There is an equivalence of 2-categories $\topstacks\simeq \hocat(\toppre)$.
\end{corollary}

\subsection{Differentiable stacks}

We will now turn to differentiable stacks.  As mentioned above, this
case is proved in precisely the same manner as the topological case.
In this case, we work with the site $\diffcat$ of small differentiable
manifolds with the \'{e}tale topology.
\begin{definition}\label{def:diff_prestack}
  A \myemph{differentiable prestack} is a prestack $A$ such that there
  exists a manifold $U$ in $\diffcat$ and a map $q\colon U\to A$
  such that, for all manifolds $V$ in $\diffcat$ and maps $f\colon V\to A$, the
  isocomma object $(f,q)$ is representable and the map $(f,q)\to V$
  is an \'{e}tale surjection.
\end{definition}
As in the topological case, we may associate to each differentiable
prestack $A$ and chart $e\colon U\to A$ a differentiable groupoid
$G^{e}$.  To such a differentiable groupoid we then have an associated
differentiable stack $Q(A,e)$ given by
\begin{align*}
  Q(A,e)_{V} & := \textnormal{\textbf{Ringed}}\bigl((\sheaves(V),C^{\infty}(V)),(\sheaves(G^{e}),C^{\infty}(U))\bigr)
\end{align*}
where the objects are morphisms of ringed toposes and the arrows are
natural isomorphisms thereof.
\begin{theorem}\label{theorem:diffpre}
  There is a system of fibrant objects on $\diffpre$ given by taking
  the local weak equivalences and with fibrant replacement given by
  the associated stack.
  \begin{proof}
    By the differentiable analogues of Lemma \ref{lemma:top_stack_equiv} and the argument given
    in the proof of Theorem \ref{theorem:toppre}, it suffices to
    construct a local weak equivalence $A\to Q(A,e)$
    for any differentiable prestack $A$ with chart
    $e\colon U\to A$.  This was done in \emph{ibid}.
  \end{proof}
\end{theorem}
\begin{corollary}
  There is an equivalence of 2-categories $\diffstacks\simeq \hocat(\diffpre)$.
\end{corollary}

\subsection{Algebraic stacks}

The case of algebraic stacks is even closer to the topological case.
In this case we work with the site $\schemes$ of small schemes
with the \'{e}tale topology.
\begin{definition}\label{def:alg_prestack}
  An \myemph{algebraic prestack} is a prestack $A$ such that the
  following conditions are satisfied:
  \begin{enumerate}
  \item The diagonal $\Delta\colon  A\to A\times A$ is representable
    and proper.
  \item There exists a scheme $U$ in $\schemes$ and a map $q\colon U\to A$
    such that, for all schemes $V$ in $\schemes$ and maps $f\colon V\to A$, the
    map $(f,q)\to V$ is an \'{e}tale surjection.
  \end{enumerate}
\end{definition}

\begin{theorem}\label{theorem:algpre}
  There is a system of fibrant objects on $\algpre$ given by taking
  the local weak equivalences and with fibrant replacement given by
  the associated stack.
  \begin{proof}
    By the algebraic analogues of Lemmas
    \ref{lemma:top_stack_rep_diag} and \ref{lemma:top_stack_equiv},
    and the argument given in the proof of Theorem
    \ref{theorem:toppre}, it suffices to construct a local weak
    equivalence $A\to Q(A,e)$ for any algebraic prestack $A$
    with chart $e\colon U\to A$. This was done in \emph{ibid}.
  \end{proof}
\end{theorem}
\begin{corollary}
  There is an equivalence of 2-categories $\algstacks\simeq \hocat(\algpre)$.
\end{corollary}

\subsection*{Acknowledgements}

We would like to thank Timothy Porter for a recent copy of his
``Menagerie'' \cite{Porter:M} and for advice regarding the literature
on stacks.  We also benefitted from discussions of this material with
Andr\'{e} Joyal.  The second author would also like to thank the AARMS, the
Department of Mathematics and Statistics at Dalhousie University, and
the Institute for Advanced Study for their support while this research
was carried out. Both authors also thank NSERC for their financial
support of this research.  The second author also received support
from NSF Grant DMS-0635607 and the Oswald Veblen Fund.

\newcommand{\SortNoop}[1]{}

\end{document}